\input amstex.tex 
\document
\magnification=\magstep1
\baselineskip 11.4pt
\NoBlackBoxes
%
%
\par \vskip 2pc \noindent
$$\text{\rm {\bf THE ALTERNATING GROUPS AND K3 SURFACES}}$$

\par \vskip 1pc \noindent
$$\text{\rm {\bf D. -Q. Zhang}}$$

\par \vskip 2pc \noindent
{\bf Abstract.}
In this note, we consider all possible extensions $G$ of a non-trivial
perfect group $H$ acting faithfully on a $K3$ surface $X$. 
The pair $(X, G)$ is proved
to be uniquely determined by $G$ if the transcendental value
of $G$ is maximum.
In particular, we have $G/H  \le ({\bold Z}/(2))^{\oplus 2}$,
if $H$ is the alternating group $A_5$ and normal in $G$. 

\par \vskip 2pc \noindent
{\bf Introduction}

\par \vskip 1pc \noindent
We work over the complex numbers field ${\bold C}$.
A {\bf K3} surface $X$ is a simply connected projective surface
with a nowhere vanishing holomorphic 2-form $\omega_X$.
In this note, we will consider finite groups in $\text{Aut}(X)$.
An element $h \in \text{Aut}(X)$ is {\bf symplectic}
if $h$ acts trivially on the 2-form $\omega_X$.
A group $G_N \subseteq \text{Aut}(X)$ is {\bf symplectic}
if every element of $G_N$ is symplectic.

\par \vskip 1pc \noindent
According to Nikulin [Ni1], Mukai [Mu1] and Xiao [Xi], there are 
exactly 80 abstract finite groups which can act symplectically
on $K3$ surfaces.
Among these 80, there are exactly four perfect groups 
($G$ is perfect if the commutator subgroup $[G, G] = G$):
$A_5, \, L_2(7), \, A_6, \, M_{20} = C_2^4 : A_5$ (the Mathieu group of degree 20),
where the first three are also the only non-abelian simple
groups which can act on a $K3$ surface symplectically,
and the last is the the symplectic finite group 
with the largest order $960$.

\par \noindent
The common thing shared by the three bigger perfect groups 
$G_N = L_2(7), A_6$ and $M_{20}$, 
is that they all can be extended
to a bigger group $G = G_N . \mu_4$ acting faithfully on
a $K3$ surface $X$. Moreover, the pair $(X, G)$ turns out to be unique
in each case, [Ko2], [OZ3], [KOZ1].

\par \noindent
So one would expect that $A_5$, being a smaller one,
should be extendable to a bigger
group $G = A_5 . \mu_I$ for some $I \ge 3$.
However, our result below shows that this is not the case.
Indeed, only $I = 1$, or $2$ is possible.

\par \vskip 1pc \noindent
{\bf Theorem A.} Suppose that a finite group $G$ acts faithfully
on a $K3$ surface. Suppose further that $G$ contains $A_5$ as a normal
subgroup. Then $G$ equals one of the following four groups,
each realizable (see Example {\bf 1.10}):
$$A_5, \,\, S_5, \,\, A_5 \times \mu_2, \,\, S_5 \times \mu_2.$$

\par \vskip 0.5pc \noindent
To be precise, as in {\bf (1.0)} below, for every
finite group $G$ acting on a $K3$ surface $X$,
the symplectic elements of $G$ (i.e., those $h$
acting trivially on the non-zero 2-form $\omega_X$)
form a normal subgroup $G_N$ such that $G/G_N \cong \mu_I$
(the cyclic group of order $I$ in ${\bold C}^*$). Namely, we have $G = G_N . \mu_I$
(see {\bf Notation} below).
The natural number $I = I(G)$ is determined by the action of $G$ on $X$ and
called the {\bf transcendental value} of (the action of) $G$.

\par \vskip 1pc \noindent
It is proved in [Ko2], [OZ3] and [KOZ1] that for the three bigger perfect groups $G_N$ above,
there is an extension $G = G_N . \mu_I$ such that the transcendental
value $I = I(G)$ equals $4$. However, for the smaller perfect (and
also simple) group $A_5$, we have:

\par \vskip 1pc \noindent
{\bf Theorem B.} Suppose that a finite group $G$ acts faithfully
on a $K3$ surface. Suppose further that $G$ contains $A_5$ as
a normal subgroup. Then the transcendental value $I(G)$ 
equals 1 or 2 (both attainable as shown in Example {\bf 1.10}).

\par \vskip 1pc \noindent
A bit more surprise comes from the next result:
the existence of action by a perfect group (together with the transcendental
value being $4$) will guarantee the existence of action by a quite large group 
$G$ as well as the uniqueness of the pair $(X, G)$.

\par \vskip 1pc \noindent
{\bf Theorem C.} Suppose that a finite group $G$ acts faithfully
on a $K3$ surface $X$. Suppose further that $G$ contains a non-trivial perfect group
$H$ as a subgroup (not necessarily normal). Then we have:

\par \noindent
(1) The transcendental value $I(G) \le 4$.

\par \noindent
(2) If $I(G) = 4$, then $G = L_2(7) . \mu_4$, $A_6 . \mu_4$ or $M_{20} . \mu_4$,
and the pair $(X, G)$ is unique, up to isomorphisms, in all three cases.

\par \vskip 1pc \noindent
{\bf Remark D.} (1) The three subgroups $L_2(7)$, $A_6$ and $M_{20}$ of $G$
in Theorem C
are all equal to $G_N$ in notation of {\bf (1.0)}, and are the only perfect
groups among the 11 maximum symplectic $K3$ groups [Mu1].
So the maximality of the transcendental value $I(G)$ in the situation of Theorem C
guarantees the maximality of the symplectic part $G_N$ of $G$.
This also shows the importance of studying non-symplectic K3 groups.

\par \noindent
(2) Regarding Theorems B and C, the readers may wonder whether the action 
of ${\widetilde A}_6 = A_6 : \mu_4$ on
a $K3$ surface $X$ induces an action of $H . \mu_4$ on $X$ with $H = A_5$ a smaller
perfect (indeed simple) group.
To elaborate, the unique group structure of $\widetilde{A}_6$ (and also the unique pair 
$(X, \widetilde{A}_6)$) is described in [KOZ1, 2]. In particular,
the natural conjugation map $\widetilde{A}_6 \rightarrow \text{Aut}(A_6)$
($x \mapsto c_x$; see {\bf Notation} below) has the Mathieu group $M_{10}$
as its image; therefore, the conjugation $\mu_4$ action switches the two different
conjugacy classes of order 3 in $A_6$ [CS, Ch 10, \S 1.5].
On the other hand, for $\widetilde{A}_6$ to contain an $A_5 . \mu_4$,
the conjugation $\mu_4$ action should stabilize at least one $A_5$ in $A_6$
and also preserve the unique conjugacy class of order 3 in this $A_5$,
which is impossible.

\par \noindent
(3) The same construction in [OZ3, Appendix] shows that there is a
smooth non-isotrivial family of $K3$ surfaces $f : {\Cal X} \rightarrow {\bold P}^1$
such that all fibres admit $A_6$ actions and infinitely many of them
are algebraic K3 surfaces. So, the symplectic part alone
can not determine the surface uniquely, and the study of transcendental value is needed.

\par \vskip 1pc \noindent
The main tools of the paper are the Lefschetz fixed point formula
(both the topological version and vector bundle version due to Atiyah-Segal-Singer 
[AS2, 3]),
the representation theory on the $K3$ lattice and the study in [Z2] 
on automorphism groups
of Niemeier lattices (in connection with Golay binary or ternary codes)
where the latter is much inspired by Conway-Sloane [CS], Kondo [Ko1] and 
Mukai [Mu2]. 

\par \noindent
The reduction to the study of automorphisms of Niemeier lattices
was pioneered by Nikulin (see e.g., [Ni3, end of section 1.14])
and further developed by Kondo (see e.g. [Ko1]).

\par \vskip 1pc \noindent
We believe that the way of combining different very powerful machinaries to 
deduce results as done in the paper should be applicable
to the study of other problems. Our humble paper also demonstrates
the powerfulness and depth of these algebraic results in the study of geometry.
The information we compute in Proposition 1.4 (and its generalization in the future) 
should be of independent interest and use in understanding the geometry of $K3$ surfaces.
  
\par \vskip 1pc \noindent
{\bf Note.} "{\bf Maple}" was used in solving the linear equations in 
the crucial Proposition 1.4. We refer to Shimada [Sh1, Sh2, Sh3]
for more computations in Algebraic Geometry.

\par \vskip 1pc \noindent
{\bf Notation.}

\par \noindent
{\bf 1.} When we write $G = G_N . \mu_I$ we mean that
$G$ acts on a $K3$ surface $X$ satisfying the situation in 
{\bf (1.0)} below.

\par \noindent
{\bf 2.} $S_n$ is the symmetric group in $n$ letters,
$A_n$ ($n \ge 3$) the alternating group
in $n$ letters and $\mu_I = \langle \text{\rm exp}(2 \pi \sqrt{-1})/I \rangle$
the multiplicative group of order $I$ in ${\bold C}^*$.

\par \noindent
{\bf 3.} For a group $G$, we write $G = A . B$ if $A$ is normal in $G$
so that $G/A = B$. We write $G = A : B$ if assume further that
$A$ is normal in $G$ and $B$ is a subgroup of $G$ such that
the composition $B \rightarrow G \rightarrow G/A = B$ is the identity
(we say then that $G$ is a {\bf semi-direct product} of $A$ and $B$).

\par \noindent
{\bf 4.} For groups $H \le G$ and $x \in G$ we denote by
$c_x : H \rightarrow G$ ($h \mapsto c_x(h) = x^{-1} h x$)
the {\bf conjugation} map.

\par \noindent
{\bf 5.} For a $K3$ surface $X$, we let $S_X$ and $T_X$ be
the Neron-Severi and transcendental lattices. 
So the $K3$ lattice $H^2(X, {\bold Z})$
contains $S_X \oplus T_X$ as a sublattice of finite index.

\par \vskip 1pc \noindent
{\bf Acknowledgement.} This work was done during
the author's visits to Hokkaido University, University of Tokyo
and Korea Institute for Advanced Study in the summer of 2004. 
The author would like to
thank the institutes and Professors I. Shimada, K. Oguiso
and J. Keum for the support and warm hospitality.

\par \vskip 2pc \noindent
{\bf \S 1. Preparations and Examples}

\par \vskip 1pc \noindent
{\bf (1.0).} In this section, we will prepare some basic results to be used
late. Let $X$ be a $K3$ surface with a non-zero 2-form $\omega_X$
and let $G \subseteq \text{Aut}(X)$ be a finite
group of automorphisms. For every $h \in G$,
we have $h^* \omega_X = \alpha(h) \omega_X$ for some scalar
$\alpha(h) \in {\bold C}^*$.
Clearly, $\alpha : G \rightarrow {\bold C}^*$
is a homomorphism. A fact in basic group theory says that $\alpha(G)$ is
a finite cyclic group, so 
$\alpha(G) = \mu_I = \langle \text{\rm exp}(2 \pi \sqrt{-1}/I) \rangle$ for some $I \ge 1$.
This natural number $I = I(G)$ is called the {\bf transcendental}
value of $G$. It is known that $I = I(G)$ for some $G$
if and only if that the Euler function $\varphi(I) \le 21$
and $I \ne 60$ [MO].

\par \noindent
Set $G_N = \text{\rm Ker}(\alpha)$. Then we have the {\bf basic exact sequence} below:
$$1 \longrightarrow G_N \longrightarrow G 
\overset{\alpha}\to{\longrightarrow} \mu_I \longrightarrow 1.$$

\par \noindent
For the $G$ in the basic exact sequence, we write $G = G_N . \mu_I$, though
there is no guarantee that $G = G_N : \mu_I$
(a semi-direct product).

\par \vskip 1pc \noindent
{\bf Fact 1.0A.} If $G$ is a finite perfect group, i.e., the commutator
group $[G, G] = G$ (especially, if $G$ is a non-abelian simple group
like $A_5$), then $G = G_N$. 

\par \noindent
{\bf 1.0B.} $G_N$ acts trivially on the transcendental lattice $T_X$
(Lefschetz theorem on $(1, 1)$-classes).

\par \noindent
{\bf 1.0C.} If a subgroup $H \le G_N$ fixes a point $P$, then
$H < SL(T_{X, P}) \cong SL_2({\bold C})$ [Mu1, {\bf (1.5)}]. 
The finite subgroups of
$SL_2({\bold C})$ are listed up in [Mu1, {\bf (1.6)}]. These are
cyclic $C_n$, binary dihedral (or quaternion) $Q_{4n}$ ($n \ge 2$), binary tetrahedral $T_{24}$, 
binary octahedral $O_{48}$ and binary icosahedral $I_{120}$.

\par \vskip 1pc \noindent
{\bf Lemma 1.1.} Suppose that $G := A_5 . \mu_I$ (with $G_N = A_5$) acts faithfully
on a $K3$ surface $X$.

\par \noindent
(1) The Picard number $\rho(X) \ge 19$, and
$I = 1, 2, 3, 4, 6$. Moreover, $\rho(X) = 20$ if $I \ge 3$.

\par \noindent
(2) We have $G = A_5 : \mu_I$, i.e.,
a semi-prodcut of a normal subgroup $A_5$ and a subgroup $\mu_I$ of $G$.
Moreover, $G = A_5 \times \mu_I$ if $I = 3$.

\par \vskip 1pc \noindent
{\it Proof.} (1) In notation of [Xi, the list], 
$\rho(X) =  \text{\rm rank} \, S_X \ge c+1 = 19$.
Also the Euler function $\varphi(I)$ divides $\text{\rm rank} \, T_X = 22 - \rho(X)$ 
by [Ni1, Theorem 0.1]. So (1) follows.

\par \noindent
(2) Let $g \in G$ such that $\alpha(g)$ is a generator of $\mu_I$.
Since $\text{Aut}(A_5) = S_5 > A_5$ and the conjugation homomorphism 
$A_5 \rightarrow \text{Aut}(A_5)$
($x \mapsto c_x$) is an isomorphism onto $A_5$, 
the conjugation map $c_g$ equals $c_{(12)a}$ or $c_{a}$ on $A_5$
for some $a \in A$. Replacing $g$ by $g a^{-1}$, we may assume that
$c_g = c_{(12)}$ or $c_{\text{\rm id}}$. Thus $g^2$ commutes with every element
in $A_5$. If $2 | I$, then $g^I \in \text{\rm Ker}(\alpha) = A_5$
is in the centre of $A_5$ (which is trivial) and hence $\text{\rm ord} (g) = I$;
thus $G = A_5 : \mu_I$. If $I = 3$, then $\text{gcd}(3, \text{\rm ord} (g)/3) = 1$
as proved in [IOZ] or [Og, Proposition 5.1]; so replacing $g$
by $g^{\ell}$ with $\ell = \text{\rm ord} (g)/3$ (or $2 \text{\rm ord} (g)/3$),
we have $G = A_5 \times \langle g \rangle = A_5 \times \mu_3$.

\par \vskip 1pc \noindent
The third result below [Ni1, \S 5] 
is crucial in classifying symplectic groups in [Mu1].
The second uses the fact $A_5 \subset \text{\rm Aut}(X)$ in an essential way.

\par \vskip 1pc \noindent
{\bf Lemma 1.2.} (1) Let $h$ be a non-symplectic involution on a $K3$ surface $X$.
Then $X^h$ is a disjoint union of $s$ smooth curves $C_i$
with $0 \le s \le 10$.
To be precise, $X^h$ (if not empty) is either a disjoint union of a genus $\ge 2$
curve $C$ and a few ${\bold P}^1$'s, or
a disjoint union of a few elliptic curves and ${\bold P}^1$'s,
or a disjoint union of a few ${\bold P}^1$'s.

\par \noindent
(2) For $h$ in (1), suppose further that $A_5
\subseteq \text{Aut}(X)$. Then $\chi_{\text{\rm top}}(X^h) \le 18$.

\par \noindent
(3) If $\delta$ is a non-trivial symplectic automorphism of finite 
order on a $K3$ surface $X$, then $\text{\rm ord} (\delta) \le 8$ and
$X^{\delta}$ is a finite set. 
To be precise, if $\text{\rm ord} (\delta) = 2, 3, 4, 5, 6, 7, 8$,
then $|X^{\delta}| = 8, 6, 4, 4, 2, 3, 2$, respectively; see [Ni1, \S 5] for the proof.
In particular, if $A_5 \subseteq \text{Aut}(X)$ then
$\sum_{\delta \in A_5} \chi_{\text{\rm top}}(X^{\delta}) = 360$ (see {\bf (1.0A)}).

\par \vskip 1pc \noindent
{\it Proof.} (1) Locally, at a point $P \in X^h$, we have
$h | P : (x, y) \rightarrow (x, -y)$ for some coordinates around $P$,
because $h$ is non-symplectic. Thus around $P$, our $X^h = \{y = 0\}$
which is smooth. For the range of $s$, see [Ni2] or [Z1].
If $X^h$ contains a genus $\ge 2$ curve $C$, then the big and nefness
of $C$ and the Hodge index theorem show that the other $s-1$
curves are negative definite, whence are ${\bold P}^1$'s. 
So (1) is true.

\par \noindent
(2) Let $X^h = \coprod_{i=1}^s C_i$ be as in (1).
Then $\chi_{\text{\rm top}}(X^h) = \sum_{i=1}^s (2 - 2g(C_i)) \le 2s \le 20$.
If (2) is false, then $s = 10$ and $C_i \cong {\bold P}^1$.
Thus, by [OZ1, Theorem 4], $X$ equals $X_4$: the unique 
$K3$ surface of Picard number $\rho(X) = 20$ and $|\text{Pic} X| = -4$.
Now $A_5 \subset \text{Aut}(X_4)$, where the latter is given
in [Vi]. This is impossible by the simplicity of $A_5$ and the
precise description of $\text{Aut}(X_4)$ there (see the proof of [KOZ1, Prop 4.1 (3)]).

\par \vskip 1pc \noindent
For an automorphism $h$ on a smooth algebraic surface $Y$, we split the
pointwise fixed locus as the disjoint union
of 1-dimensional part and the isolated part:
$Y^h = Y^h_{1-\text{\rm dim}} \coprod Y^h_{\text{\rm isol}}$.
The proof of (1) below is similar to that for (1) in {\bf (1.2)}.

\par \vskip 1pc \noindent
{\bf Fact 1.3.} (1) $Y^h_{1-\text{\rm dim}}$ (if not empty) is a disjoint union of smooth curves.

\par \noindent
(2) The Euler number $\chi_{\text{\rm top}}(Y^h_{1-\text{\rm dim}}) = \sum_{C}
(2 - 2 g(C)) = 2n_h$ for some integer $n_h$, where
$C$ runs in the set $Y^h_{1-\text{\rm dim}}$ of curves.

\par \noindent
(3) The Euler number $\chi_{\text{\rm top}}(Y^h) = m_h + 2 n_h$, where
$m_h = |Y^h_{\text{\rm isol}}|$.

\par \vskip 1pc \noindent
The results of [IOZ] below follow from the application of
Lefschetz fixed point 
formula to the trivial vector bundle in Atiyah-Segal-Singer
[AS2, AS3, pages 542 and 567]. The results themselves should
be very useful and informative for other studies in the future.

\par \vskip 1pc \noindent
{\bf Important Proposition 1.4.} Let $X$ be a $K3$ surface and let
$h \in \text{Aut}(X)$ be of order $I$ such that $h^* \omega_X = \eta_I \omega_X$
for some primitive $I$-th root $\eta_I$ of 1.

\par \noindent
(1) Suppose that $I = 3$. Then $m_h = 3 + n_h$ and hence 
$\chi_{\text{\rm top}}(X^h) = 3(1 + n_h)$. Moreover, $-3 \le n_h \le 6$.

\par \noindent
(2) Suppose that $I = 4$. Then $m_h = 4 + 2 n_h$ and hence 
$\chi_{\text{\rm top}}(X^h) = 4(1 + n_h)$. Moreover, $-2 \le n_h \le 4$.

\par \noindent
(3) Suppose that $I = 3$, or $4$. 
If $\delta \in \text{Aut}(X)$ is symplectic of order 5 and
commutes with $h$. Then $|X^{h \delta}| = 4$.

\par \noindent
(4) Suppose that $I = 4$. If $\delta \in \text{Aut}(X)$ is symplectic of order 3 and
commutes with $h$ then $6 \ge |X^{h^2 \delta}| \ge |X^{h \delta}|
\in \{2, 4, 6\}$.

\par \vskip 1pc \noindent
{\it Proof.} (1) The first part is proved in [OZ1, Lemma 2.3].
Note that $h^* | T_X$ can be diagonalized as 
$\text{\rm diag}[\eta_3, \eta_3^2]^{\oplus s}$ ($s \ge 1$) by [Ni1, Theorem 0.1].
So as in {\bf (1.7)} below,
$\chi_{\text{\rm top}}(X^h) = 2 + \text{\rm Tr} (h^*|T_X) + \text{\rm Tr}(h^* | S_X)
\le 2 -s + \text{\rm rank} \, S_X \le 21$, whence $n_h \le 6$. 
Also $m_h \ge 0$ implies that $n_h \ge -3$.

\par \noindent
(2) As in [OZ1, Lemma 2.3], we calculate the holomorphic Lefschetz number $L(h)$ in two
ways as in [AS2, 3, pages 542 and 567],
where $X^h_{\text{\rm isol}} = \{P_j | 1 \le j \le m_h\}$ 
(so $h^*|T_{P_j} = (\eta_4^{-1}, \eta_4^2)$
up to switching the coordinates of the tangent plane at $P_j$),
$X^h_{\text{\rm 1-dim}} = \{C_k\}$,
$gC_k = g(C_k)$ the genus, and $\eta_4^{-1}$ the eigenvalue
of the action $h_*$ on the normal bundle of $C_k$
(in the first equation below
we used Serre duality, while the last is from the first two with $x = \eta_4$):
$$\gather
L(h) = \sum_{i=0}^2 (-1)^i \text{\rm Tr} (h^* | H^i(X, {\Cal O}_X)) = 1 + \eta_4^{-1}, \\
L(h) = \sum_{j=1}^{m_h} a(P_j) + \sum_k b(C_k), \\
a(P_j) = 1/\text{\rm det}(1 - h^*|T_{P_j}) = 1/(1 - \eta_4^{-1})(1 - \eta_4^2), \\
b(C_k) = (1 - gC_k)/(1 - \eta_4) - \eta_4 C_k^2/(1 - \eta_4)^2 = (1- gC_k)(1 + \eta_4)/(1 - \eta_4)^2, \\
0 = -(1 + x^{-1}) + m_h/(1 - x^{-1})(1 - x^2) + n_h(1+x)/(1-x)^2.
\endgather$$

\par \noindent
Noting that $x = \eta_4$ satisfies $x^2 = -1$ and solving the last equation,
we get $m_h = 4 + 2 n_h$. The second part of (2) is similar to (1),
noting that $h^* | T_X$ can be diagonalized as $\text{\rm diag}[\eta_4, - \eta_4]^{\oplus s}$ 
($s \ge 1$).

\par \noindent
(3) $\&$ (4). In (4), note that $X^{h^i \delta} = X^{h^i} \cap X^{\delta}$ ($i = 1, 2$).
So the inequalities there hold and we have only to calculate $|X^{h \delta}|$;
see {\bf (1.2)}.

\par \noindent
Let $g \in \text{\rm Aut}(X)$ such that $\text{\rm ord}(g) = k I$ and 
$g^* \omega_X = \eta^k \omega_X$
where $\eta = \eta_{kI}$ is a primitive $kI$-th root of 1. 
(We set $g = h \delta$ in (3) and (4).) If $k \ge 2$ and $\text{\rm gcd}(k, I) = 1$ (these are true in (3) and (4)),
then $g^I$ is of order $k$ and symplectic, so
$X^g \subseteq X^{g^I}$ is a finite set by {\bf (1.2)}. Namely, $X^g = X^g_{\text{\rm isol}} 
= \{P_j | 1 \le j \le m_g\}$ say. Let $M_g(i)$ be the set
of points $P$ in $X^g$ satisfying $g^*|T_P = (\eta^{-i}, \eta^{k + i})$
(up to switching the coordinates of the tangent plane at $P$;
so $a(P) = 1/(1 - \eta^{-i})(1 - \eta^{k + i})$ in the notation for the formula of $L(g)$).
Put $m_g(i) = |M_g(i)|$.
Then for $(I, k) = (3, 5)$ (the first case in (3)),
we have $X^g = \coprod M_g(i)$ and $m_g = \sum_i m_g(i)$, where $i \in \{1, \dots, 4, 11, 12\}$;
for $(I, k) = (4, 5)$ (the second case in (3)), we have 
$m_g = \sum_i m_g(i)$, where $i \in \{1, \dots, 4, 6, 7, 16, 17\}$;
for $(I, k) = (4, 3)$ (the case in (4)), we have 
$m_g = \sum_i m_g(i)$, where $i \in \{1, 2, 4, 10\}$.

\par \noindent
As in (2), we have the following, where $x = \eta = \eta_{kI}$ and
$i$ runs in the set specified above:
$$(*) \hskip 1pc 0 = -(1 + x^{-k}) + \sum_i \sum_{P \in M_g(i)} a(P) =
-(1 + x^{-k}) + \sum_i m_g(i) /(1 - x^{-i})(1 - x^{k + i}).$$

\par \noindent
For $(I, k) = (3, 5)$, $x$ satisfies the minimal polynomial 
$\Phi_g(x) = 1 - x + x^3 - x^4 + x^5 - x^7 + x^8$
and also $x^{15} = 1$, $x^{10} = -1 - x^5$.
Substituting these into the equation (*) multiplied by the common denomenator
(which is not zero), we will get an equation of degree $\le 7$ in $x$
with coefficients linear in $m_g(i)$. The minimality of $\Phi_g(x)$ implies that
all 8 coefficients are zero. Solving these 8 linear equations,
we obtain, where $m_i = m_g(i)$:
$$(**) \hskip 1pc m_1 = m_4, \,\, m_2 = -1 + m_3, \,\, m_{11} = -1 + m_4, \,\,m_{12} = m_3.$$

\par \noindent
By {\bf (1.2)}, we have
$4 = m_{g^3} \ge m_g = \sum_{i=1}^4 m_i + \sum_{i=11}^{12} m_i = -2 + 3(m_3+m_4)$.
So $m_3 + m_4 \le 2$. This together with the condition $m_i \ge 0$ and the
relations in (**), imply that 
$[m_1, m_2, m_3, m_4, m_{11}, m_{12}] = [1, 0, 1, 1, 0, 1]$. 
In particular, $m_g = 4$.

\par \noindent
For $(I, k) = (4, 5)$, $x$ satisfies the minimal polynomial 
$\Phi_g(x) = 1 - x^2 + x^4 - x^6 + x^8$
and also $x^{20} = 1$, $x^{10} = -1$. As above, solving $(*)$, 
we obtain, where $m_i = m_g(i)$:
$$
\gather (***)
\hskip 0.5pc m_1 = -3 + 2m_3 -3m_4 + 4m_6 - 2m_7, \,\,
m_2 = -1 + m_3 - 2m_4 + 2m_6, \\
m_{16} = -5 + 2m_3 - 4m_4 + 5m_6 - 2m_7, \,\,
m_{17} = 3 + 2m_4 - 2m_6 + m_7.
\endgather$$

\par \noindent
One can check that the following is the only possibility of $m_i$ 
satisfying the relations in $(***)$ and that $0 \le m_i \le m_g
\le m_{g^4} = 4$ by {\bf (1.2)}; in particular, $m_g = 4$:
$$[m_1, m_2, m_3, m_4, m_6, m_7, m_{16}, m_{17}]
= [1, 1, 0, 0, 1, 0, 0, 1].$$

\par \noindent
For $(I, k) = (4, 3)$, $x$ satisfies the minimal polynomial 
$\Phi_g(x) = 1 - x^2 + x^4$
and also $x^{12} = 1$, $x^{6} = -1$. As above, solving $(*)$, 
we obtain, where $m_i = m_g(i)$:
$$
(****) \hskip 1pc m_1 = 3 + 3m_2 - 2m_4, \,\, m_{10} = 1 + 2m_2 - m_4.
$$

\par \noindent
One can check that the following are the only possibilities of $m_i$ 
satisfying the relations in $(****)$ and $0 \le m_i \le m_g
\le m_{g^4} = 6$, {\bf (1.2)}; in particular, $m_g = 2, 4, 6$
(so {\bf (1.4)} is done):
$$[m_1, m_2, m_4, m_{10}]
= [3, 0, 0, 1], \,\, [1, 0, 1, 0], \,\, [2, 1, 2, 1], \,\, [0, 1, 3, 0].$$

\par \vskip 1pc \noindent
The following two results can be found in [Ni1, Theorem 0.1], [MO, Lemma (1.1)],
or [OZ3, Lemma (2.8)].

\par \vskip 1pc \noindent
{\bf Lemma 1.5.} Suppose that $X$ is a $K3$ surface of Picard
number $\rho(X) = 20$ and $g$ an order-4 automorphism
such that $g^* \omega_X = \eta_4 \omega_X$ with a primitive
4-th root $\eta_4$ of 1. Then we can express the transcendental lattice
$T_X$ as $T_X = {\bold Z}[t_1, t_2]$ so that
$t_2 = g^*(t_1)$, $g^*(t_2) = -t_1$. In particular,
the intersection form $(t_i . t_j) = \text{\rm diag}[2m, 2m]$ for some
$m \ge 1$.

\par \vskip 1pc \noindent
Now we assume that $G = G_N . \mu_I$ (with $I = I(G)$)
acts on a $K3$ surface $X$.
When $G_N = A_5$,
we will determine the action of $G_N$
on the Neron Severi lattice $S_X$ of $X$:

\par \vskip 1pc \noindent
{\bf Lemma 1.6.} (1) Suppose that $A_5$ acts on a $K3$ surface $X$, 
and $\text{\rm rank} \, S_X = 20$ 
(this is true if $I \ge 3$
by {\bf (1.1)}). Then
we have the irreducible decomposition below (in the
notation of Atlas for the characters of $A_5$),
where $\chi_1$ (the trivial character), 
$\chi_4$ and $\chi_5$ have dimensions 1, 4 and 5, respectively,
where $\chi_i'$ is a copy of $\chi_i$:
$$S_X \otimes {\bold C}
= \chi_1 \oplus \chi_1' \oplus \chi_4 \oplus \chi_4' \oplus \chi_5 \oplus \chi_5'.$$

\par \noindent
(2) For conjugacy class $nA$ (and $nB$) of order $n$ in $A_5$ and the characters 
$\chi_i$ of $A_5$, we have the following {\bf Table 1} from [Atlas], where $Z$ is respectively 
$1A$, $2A$, $3A$, $5A$ or $5B$:
$$\gather
[\chi_1, \chi_2, \chi_3, \chi_4, \chi_5](Z) =
[1, 3, 3, 4, 5], \hskip 1pc [1, -1, -1, 0, 1], \hskip 1pc [1, 0, 0, 1, -1], \\
[1, \, (1- \sqrt{5})/2, \,\,(1 + \sqrt{5})/2, \, -1, 0], \hskip 1pc
[1, \, (1+ \sqrt{5})/2, \,\,(1 - \sqrt{5})/2, \, -1, 0].
\endgather $$

\par \vskip 1pc \noindent
{\it Proof.} Applying the Lefschetz fixed point formula to the action
of $A_5$ on 
$H^*(X, {\bold Z}) = \oplus_{i=0}^4 H^i(X, {\bold Z})$
and noting that $H^2(X, {\bold Z})$
contains $S_X \oplus T_X$ as a finite index sublattice, we obtain
(see also {\bf (1.0A-B)} and {\bf (1.2)}):
$$2 + \text{\rm rank} \, T_X + \text{\rm rank} (S_X)^{A_5} = 
\text{\rm rank} \,H^*(X, {\bold Z})^{A_5} = \frac{1}{|A_5|} \sum_{a \in A_5} \chi_{\text{\rm top}}(X^a) = 360/60 = 6.$$

\par \noindent
Thus $\text{\rm rank} \, S_X^{A_5} = 2$.
So the irreducible decomposition is of the
following form, where $a_{i}$ are non-negative integers:
$$S(X) \otimes {\bold C} = 2 \chi_{1} \oplus a_{2}\chi_{2} 
\oplus a_{3}\chi_{3} \oplus a_{4}\chi_{4} \oplus a_{5}\chi_{5}.$$
As in {\bf (1.7)} below, using the topological Lefschetz
fixed point formula, the fact that $\text{\rm rank}\, T(X)$
$= 2$ and {\bf (1.0B)}, we have, for $a \in A_{5}$, that:
$$\chi_{\text{\rm top}}(X^{a}) = 4 + \text{\rm Tr}(a^{*} \vert S(X))$$
Running $a$ through the five conjugacy classes and
calculating both sides, using {\bf (1.2)} and the character Table 1 in (2),
we obtain the following system of equations:
$$\gather 20 = 2 + 3(a_{2} + a_{3}) + 4a_{4} + 5a_{5}, \\
4 = 2 - (a_{2} + a_{3}) + a_{5}, \\
2 = 2 +  a_4 - a_{5}, \\
0 = 2 + \frac{1 - \sqrt{5}}{2}a_{2} + \frac{1 + \sqrt{5}}{2} a_{3} 
- a_{4}, \\
0 = 2 +  \frac{1 + \sqrt{5}}{2}a_{2} + \frac{1 - \sqrt{5}}{2} a_{3} 
- a_{4}.
\endgather$$

\par \noindent
Now, we get the result by solving this system of Diophantine equations.

\par \vskip 1pc \noindent
{\bf (1.7).} Note that $\text{Aut}(A_5) = S_5$. For a group $G = A_5 . \mu_I$ 
(and the map $\alpha$) in {\bf (1.0)},
we have the natural homomorphism below, which is injective
(since its restriction on $A_5$ is an injection onto $A_5 \times (1)$),
where $c_x : a \mapsto c_x(a) = x^{-1} a x$ is the conjugation map:
$$\gather
G \longrightarrow \text{Aut}(A_5) \times \mu_I = S_5 \times \mu_I, \\
x \mapsto (c_x, \alpha(x)).
\endgather$$

\par \vskip 0.5pc \noindent
{\bf Lemma.} Suppose that $G = A_5 . \mu_4$ acts on a $K3$ surface $X$
(i.e., $G_N = A_5$ and $I(G) = 4$). 
Then $G = A_5 : \mu_4$, but $G \ne A_5 \times \mu_4$. 
Our $G \rightarrow S_5 \times \mu_4$
($(x \mapsto (c_x, \alpha(x))$) is an injective homomorphism
and the group structure of such $G$ is unique up to isomorphisms.

\par \vskip 1pc \noindent
{\it Proof.} By {\bf (1.1)}, we have $G = A_5 : \mu_4$.
Suppose the contrary $G = A_5 \times \mu_4$.
Write $\mu_4 =$ $\langle g \rangle$. In notation of {\bf (1.6)},
the $g$ either stabilizes $\chi_i$ or swtiches $\chi_i$ with $\chi_i'$ ($i = 4$ or $5$;
then denoted as $\chi_i \overset{g}\to{\leftarrow \rightarrow} \chi_i$,
and $\text{\rm Tr}(g^* | (\chi_i \oplus \chi_i')) = 0$)).
Since $G$ stabilizes an ample line bundle (the pull back of 
an ample line bundle on $X/G$) and since $G$ acts on $S_X^{A_5}$
(whose ${\bold C}$-extension is $\chi_1 \oplus \chi_1'$), we may assume that
$g^*|(\chi_1 \oplus \chi_1') = \text{\rm diag}[1, \pm 1]$ w.r.t. to
a suitable basis.
If $\chi_i$ is $g$-stable then $g^* | \chi_i$ is a scalar $\zeta_4^c$
with $\zeta_4 = \text{\rm exp}(2 \pi \sqrt{-1}/4)$,
by Schur's lemma.

\par \noindent
Let $a \in A_5$. Then $(g a)^* | T_X = g^* | T_X$ (see {\bf (1.0B)}) and the latter
can be diagonalized as $\text{\rm diag}[\zeta_4, \zeta_4^{-1}]$
by [Ni1, Theorem 0.1] and {\bf (1.1)}.  Hence $\text{\rm Tr}(ga)^* | T_X = 0$.
By the topological Lefschetz fixed point formula
and noting that $H^2(X, {\bold Z})$ contains $S_X \oplus T_X$ as a sublattice
of finite index, we have 
$\chi_{\text{\rm top}}(X^{ga}) = \oplus_{i=0}^4 \text{\rm Tr}(ga)^* | H^i(X, {\bold Z})$
$= 2 + \text{\rm Tr}  (ga)^* | S_X + \text{\rm Tr}  (ga)^* | T_X = 2 + \text{\rm Tr}  (ga)^* | S_X$.
For $a = 5A$ (an order-5 element) in $A_5$, by {\bf (1.4)} and Table 1
in {\bf (1.6)}) (and Schur's lemma), we have:
$4 = \chi_{\text{\rm top}}(X^{g 5A}) = 2 + \text{\rm Tr}(g^* | \chi_1 \oplus \chi_1') 
+ \text{\rm Tr}  (g 5A)^* | (\chi_4 \oplus \chi_4) + 0$,
so one of the following cases occurs (using Schur's lemma):

\par \vskip 1pc \noindent
Case(i). $g^* | S_X \otimes {\bold C} = \text{\rm diag}[1, -1, -I_4, -I_4, ?, ?]$,

\par \noindent
Case(ii). $g^* | S_X \otimes {\bold C} = \text{\rm diag}[1, 1,\chi_4 \overset{g}\to{\leftarrow \rightarrow} \chi_4, ?, ?]$,

\par \noindent
Case(iii). $g^* | S_X \otimes {\bold C} = \text{\rm diag}[1, 1, I_4, -I_4, ?, ?]$,

\par \noindent
Case(iv). $g^* | S_X \otimes {\bold C} = 
\text{\rm diag}[1, 1, \zeta_4 I_4, \zeta_4^{-1} I_4, ?, ?]$.

\par \vskip 1pc \noindent
By {\bf (1.4)}, we have $(*) : \,\, 
-4 \le \chi_{\text{\rm top}}(X^g) = 4(1+n_g) = 0$ (mod 4) with 
$-2 \le n_g \le 4$. So 
$\chi_{\text{\rm top}}(X^g) =  4$ in Cases (ii), (iii) and (iv)
(using Schur's lemma). 
Thus $n_g = 0$ and $m_g = 4 + 2 n_g = 4$ by {\bf (1.4)}. 
Now $A_5$ (commuting with $g$) acts on the four isolated points $P_i$ in $X^g$,
whence fixing these four points (see {\bf (1.8)} below). So 
$A_5 < SL(T_{X, P_1})$, contradicting {\bf (1.0C)}.
In Case(i), by the fact (*) above and Schur's lemma, we have 
$\chi_{\text{\rm top}}(X^g) = 2 + (1 - 1 - 4 - 4 + 5 + 5) = 4$, 
which will lead to the same contradiction.

\par \noindent
By the proof of {\bf (1.1)} and the result in the above
paragraph, we may assume that there is an order-4 element $g \in G$
such that $\alpha(g)$ is the generator of $\mu_4$, so that
$G = A_5 : \langle g \rangle = A_5 : \mu_4$ and
the conjugation map $c_g = c_{(12)}$ on $A_5$.
Clearly, the group structure of $G$ is unique.
The lemma is proved.

\par \vskip 1pc \noindent
The two results below are either easy or well known and will be
frequently used in the arguments of the subsequent sections.

\par \vskip 1pc \noindent
{\bf Lemma 1.8.} Let $f : A_5 \rightarrow S_r$ ($r \ge 2$) be a
homomorphism.

\par \noindent
(1) If $r = 2$, $3$, or $4$, then $f$ is trivial.

\par \noindent
(2) If $\text{\rm Im}(f)$ is a transitive subgroup of the full symmetry group $S_r$
in $r$ letters $\{1, 2, \dots, r\}$ (whence $r \ge 5$ by (1)), then
$r | |A_5|$ with $|A_5|/r$ equal to the order of the subgroup of 
$A_5$ stabilizing the letter 1, so $r \in \{5, 6, 10, 12, 15, 20, 30\}$. 

\par \vskip 1pc \noindent
{\bf Lemma 1.9.}  
(1) $\text{Aut}({\bold P}^1)$ includes $A_5$ but not $S_5$ [Su, Theorem 6.17].

\par \noindent
(2) If $\text{\rm id} \ne f \in \text{Aut}({\bold P}^1)$ is an automorphism of finite order,
then $f$ fixes exactly two point of ${\bold P}^1$
(by the diagonalization of a lifting of $f$ to $GL_2({\bold C})$).

\par \noindent
(3) If $f_r$ ($r = 2$ or $3$) is an order$-r$ automorphism of an 
elliptic curve $E$, then either $f_r$ acts freely on $E$, or
the fix locus satisfies
$|X^{f_r}| = 4$ (resp. $= 3$) if $r = 2$ (resp. $r = 3$) (by the Hurwitz formula).

\par \vskip 1pc \noindent
The examples below are to show the existence of the groups in Theorems A and B.

\par \vskip 1pc \noindent
{\bf Example 1.10.} (1) {\bf $G = G_N  . \mu_I = S_5 \times \mu_2$ 
(with $G_N = S_5$, $I = 2$) acts on a $K3$.}

\par \noindent
Let $X = \{ \sum_{i=1}^5 X_i = \sum_{i=1}^6 X_i^2 = \sum_{i=1}^5 X_i^3 = 0\}
\subset {\bold P}^5$. We define the symplectic action of $\sigma \in S_5$ on $X$
(the same as in [Mu1, $n^{\circ} 3$]) and a non-symplectic involution
$g$ on $X$ as follows (see [Mu1, Lemma 2.1]):
$$\gather
\sigma : [X_1 : \dots : X_6] \mapsto [X_{\sigma(1)} : \dots : X_{\sigma(5)} :
(\text{\rm sign} \, \sigma) X_6], \\
g : [X_1 : \dots : X_6] \mapsto [X_1 : \dots : X_5 :
- X_6].
\endgather$$
Let $G = \langle S_5, g \rangle$. Then $G = S_5 \times \langle g \rangle$
is the required one.

\par \vskip 1pc \noindent
(2) {\bf $G = G_N  . \mu_I = A_5 : \mu_2 = S_5$ 
(with $G_N = A_5$, $I = 2$) acts on a $K3$ surface.}

\par \noindent
Let $X = \{ \sum_{i=1}^6 X_i = \sum_{i=1}^6 X_i^2 = \sum_{i=1}^6 X_i^3 = 0\}
\subset {\bold P}^5$. We define the action of $\sigma \in S_6$ on $X$
(the same as in [Mu1, $n^{\circ} 2$]):
$$\sigma : [X_1 : \dots : X_6] \mapsto [X_{\sigma(1)} : \dots : X_{\sigma(6)}].$$
Since $A_6$ is perfect, its action on $X$ is symplectic {\bf (1.0A)}.
If we let $\widetilde{G} = S_6$, then $\widetilde{G} = \widetilde{G}_N . \mu_2$
with $\widetilde{G}_N = A_6$ and $I = 2$ (see [Mu1, Lemma 2.1]).
Now a subgroup $G = S_5$ of $\widetilde{G}$ is the required one.

\par \vskip 2pc \noindent
{\bf \S 2. The determination of some topological invariants}

\par \vskip 1pc \noindent
Let $X$ be a $K3$ surface with a faithful action by a group of the form
$G := A_5 . \mu_4$ as in {\bf (1.0)}.
Then $G = A_5 : \mu_4$
and the {\it unique} group structure of such $G$ is given in {\bf (1.7)}.

\par \noindent
We will use the notation in {\bf (1.6)}.
Let $g$ be a generator of $\mu_4 < G$.
We may also assume the following is true (after a change of $g$):

\par \vskip 1pc \noindent
{\bf Lemma 2.1.} (1) The conjugation action $c_g(.) = c_{(12)}(.)$ on $A_5$.
So $\langle g^2 \rangle$ is in the centre of $G$ and
$G \rightarrow \text{Aut}(A_5) = S_5$ ($x \mapsto c_x$) induces an isomorphism
$G/\langle g^2 \rangle \cong S_5$.

\par \noindent
(2) $g^* \omega_X = \zeta_4 \omega_X$ with 
$\zeta_4 = \text{\rm exp}(2 \pi \sqrt{-1}/4)$.

\par \noindent
(3) $g^2$ is a non-symplectic involution on $X$
and commutes with every element in $A_5$.

\par \noindent
(4) Set $\sigma = (12)(34)$ and $\tau = (345)$. Then
$g$ commutes with every element in $\langle \sigma, \tau \rangle = S_3$.
So $G = A_5 : \mu_4 > S_3 \times \mu_4$.

\par \noindent
(5) Set $\sigma = (12)(34)$, $\gamma = (123)$. Then $g$ normalizes
$\langle \sigma, \gamma \rangle = A_4$. So $G = A_5 : \mu_4 > A_4 : \mu_4$.
Set $\sigma_1 = \sigma$
and $\sigma_2 = (13)(24)$ (all in $A_4$).

\par \noindent
(6) $g$ stabilizes both $\chi_1$ and $\chi_1'$;
the restrictions $g^* | \chi_1 = \text{\rm id}$
and $g^* | \chi_1' = \pm \text{\rm id}$ (after a change of basis).

\par \noindent
(7) $g$ either stabilizes both $\chi_4$ and $\chi_4'$
(so the restrictions of $g^*$ on $\chi_4$ and $\chi_4'$
are scalar multiplications),
or switches $\chi_4$ with $\chi_4'$.

\par \noindent
(8) $g$ either stabilizes both $\chi_5$ and $\chi_5'$
(so the restrictions of $g^*$ on $\chi_5$ and $\chi_5'$
are scalar multiplications),
or switches $\chi_5$ with $\chi_5'$.

\par \noindent
(9) Both $g^2 | \chi_i$ and $g^2 | \chi_i'$ ($i = 4, 5$) are scalar 
multiplications.

\par \vskip 1pc \noindent
{\it Proof.} (1) is from the last part of the proof of {\bf (1.7)}.
The (2) is true because $g$ is a generator of $\mu_4 < G = A_5 : \mu_4$.
The (3), (4) and (5) follow from (1). The (6) is true because
$G = A_5 : \langle g \rangle$ stabilizes one ample line bundle
(the pull back of an ample line bundle on $X/G$) and
$g$ acts on $S_X^{A_5}$ (defined over ${\bold Z}$)
whose ${\bold C}$-extension is $\chi_1 \oplus \chi_1'$.
(7), (8) and (9) are from the form of the decomposition in {\bf (1.6)}
and Schur's lemma.

\par \vskip 1pc \noindent
In the rest of the section, we will prove
the Key result {\bf (2.2)} below which will be used in
the proof of Theorems A, B and C in \S 5 and
is the consequence of {\bf (2.6)}-{\bf (2.9)} below.
The representation theory (mainly on $A_5$) is fully applied.
We divide into cases according to whether $g$ stablilizes or switches
$\chi_i$ ($i = 4, 5$).

\par \vskip 1pc \noindent
{\bf Key Proposition 2.2.} Suppose that $G = A_5 : \mu_4$ 
acts on a $K3$ surface $X$.
Then with the notation in {\bf (2.1)} and {\bf (1.4)},
$(n_g, m_g; \,\, \chi_{\text{\rm top}}(X^g), \chi_{\text{\rm top}}
(X^{g \tau}), \chi_{\text{\rm top}}(X^{g^2 \tau})$,
$\chi_{\text{\rm top}}(X^{g^2}))$ is one of the following:
$$(1, 6; 8, 2, 6, 0), \hskip 0.5pc (0, 4; 4, 4, 6, 0), 
\hskip 0.5pc (-1, 2; 0, 6, 6, 0).$$

\par \vskip 1pc \noindent
The result below is used in {\bf (2.4)} to determine the representation
of $S_3 \times \mu_4 < G$ there.

\par \vskip 1pc \noindent
{\bf Lemma 2.3.} (1) Suppose that $g$ stabilizes $\chi_4$.
Then w.r.t. to one and the same basis $\{v_1, \dots, v_4\}$, we have the 
following matrix representation of $A_4 : \mu_4$ on $\chi_4$:
$$\sigma_1^* = \text{\rm diag} [1, 1, -1, -1], \hskip 1pc
\sigma_2^* = [1, -1, 1, -1],$$
$$\gamma^* = \pmatrix 
1 & 0 & 0 & 0 \\
0 & 0 & 0 & \beta_4 \\
0 & \beta_2 & 0 & 0 \\
0 & 0 & \beta_3 & 0 \endpmatrix, \hskip 1pc
g^* = \pmatrix 
\alpha_1 & 0 & 0 & 0 \\
0 & \alpha_2 & 0 & 0 \\
0 & 0 & 0 & \alpha_5 \\
0 & 0 & \alpha_4 & 0 \endpmatrix.$$

\par \noindent
We have exactly the same kind of matrix representation of $A_4 : \mu_4$
w.r.t.
one and the same basis $\{v_1', \dots, v_4'\}$ of $\chi_4'$. But
we use $\beta_i'$ and $\alpha_i'$ for $\gamma^* | \chi_4'$ and
$g^* | \chi_4'$ instead.

\par \noindent
(2) Suppose $g$ stabilizes $\chi_5$.
Then w.r.t. to one and the same basis $\{y_1, \dots, y_5\}$, we have the 
following matrix representation of $A_4 : \mu_4$ on $\chi_5$,
where $\eta_3$ is a primitive 3rd root of 1:
$$\sigma_1^* = \text{\rm diag} [1, 1, 1, -1, -1], \hskip 1pc
\sigma_2^* = [1, 1, -1, 1, -1],$$
$$\gamma^* = \pmatrix 
\eta_3 & 0 & 0 & 0 & 0 \\
0 & \eta_3^2 & 0 & 0 & 0 \\
0 & 0 & 0 & 0 & b_5 \\
0 & 0 & b_3 & 0 & 0 \\
0 & 0 & 0 & b_4 & 0
\endpmatrix, \hskip 1pc
g^* = \pmatrix 
0 & a_2 & 0 & 0 & 0 \\
a_1 & 0 & 0 & 0 & 0 \\
0 & 0 & a_3 & 0 & 0 \\
0 & 0 & 0 & 0 & a_5 \\
0 & 0 & 0 & a_4 & 0
\endpmatrix.$$

\par \noindent
We have exactly the same kind of matrix representation of $A_4 : \mu_4$ 
w.r.t.
one and the same basis $\{y_1', \dots, y_5'\}$ of $\chi_5'$. But
we use $b_i'$ and $a_i'$ for $\gamma^* | \chi_5'$ and $g^* | \chi_5'$ instead.

\par \vskip 1pc \noindent
{\it Proof.} This follows from the character Table 1 in {\bf (1.6)} and
the fact that
the conjugation $c_g$ fixes $\sigma_1$, and
exchanges $\sigma_2$ with $\sigma_1 \sigma_2$
and $\gamma$ with $\gamma^{-1}$.

\par \vskip 1pc \noindent
{\bf Lemma 2.4.} (1) Suppose that $g$ stabilizes $\chi_4$.
Then w.r.t. to one and the same basis $\{u_1, \dots, u_4\}$, we have the 
following matrix representation of $S_3 \times \mu_4$ on $\chi_4$,
where $\eta_3$ is a primitive 3rd root of 1.
Moreover, $d_3 = \pm d_1$ and $(g^2)^* | \chi_4 = d_1^2 \,\, \text{\rm id}$:
$$\gather
\tau^* = [1, 1, \eta_3, \eta_3^2], \hskip 1pc g^* = \text{\rm diag}[d_1, -d_3, d_3, d_3],\\
\sigma^* = \text{\rm diag}[1, -1, \pmatrix 0 & 1 \\ 1 & 0 
\endpmatrix].
\endgather$$

\par \noindent
We have exactly the same kind of matrix representation of $S_3 \times \mu_4$ w.r.t.
one and the same basis $\{u_1', \dots, u_4'\}$ of $\chi_4'$. But
we use $d_i'$ for $g^* | \chi_4'$ instead.

\par \noindent
(2) Suppose that $g$ stabilizes $\chi_5$.
Then w.r.t. to one and the same basis $\{x_1, \dots, x_5\}$, we have the 
following matrix representation of $S_3 \times \mu_4$ on $\chi_5$,
where $\eta_3$ is a primitive 3rd root of 1.
Moreover, $e_2 = \pm e_1$, $(g^2)^* | \chi_5 = e_1^2 \, \text{\rm id}$
(and $e_1$ equals $a_3$ in {\bf (2.3)}):
$$\gather
\tau^* = \text{\rm diag}[1, \eta_3, \eta_3^2, \eta_3, \eta_3^2], \hskip 1pc
g^* = [e_1, e_2, e_2, -e_2, -e_2],\\
\sigma^* = \text{\rm diag}[1, \pmatrix 0 & 1 \\ 1 & 0 \endpmatrix, 
\pmatrix 0 & 1 \\ 1 & 0 \endpmatrix].
\endgather$$

\par \noindent
We have exactly the same kind of matrix representation of $S_3 \times \mu_4$ 
w.r.t.
one and the same basis $\{x_1', \dots, x_5'\}$ of $\chi_5'$. But
we use $e_i'$ for $g^* | \chi_5'$, instead.

\par \vskip 1pc \noindent
{\it Proof.} (1) follows from the character Table 1 in
{\bf (1.6)} and the fact that
$g$ commutes with both $\sigma, \tau$, if we claim only
$g^* | \chi_4 = \text{\rm diag} [d_1, d_2, d_3, d_3]$ instead.
It suffices to show that $d_2 = -d_3$.
On the one hand, over the eigenspace $V_4(\sigma = -1) \subset
\chi_4$ of $\sigma$ corresponding to the eigenvalue $-1$,
we have $g^* | V_4(\sigma=-1) = \text{\rm diag}[d_2, d_3]$.
On the other hand, by {\bf (2.3)}, 
$g^* | V_4(\sigma=-1) = \text{\rm diag}[\sqrt{\alpha_4 \alpha_5},
-\sqrt{\alpha_4 \alpha_5}]$. Thus $d_2 = -d_3$. Now $d_1 = \pm d_3$ follows
from the fact that $(g^2)^* | \chi_i$ is a scalar.

\par \noindent
(2) is similar, except the determination of $e_i$ in
$g^* = \text{\rm diag}[e_1, e_2, e_2, e_4, e_4]$.
Indeed, comparing the diagonalization in {\bf (2.3)} and here
we see also that $\text{\rm diag}[e_2, e_4] = g^* | V_5(\sigma=-1) = 
\text{\rm diag}[\sqrt{a_4a_5}, -\sqrt{a_4a_5}]$, whence $e_4 = -e_2$. 
Taking trace in {\bf (2.3)} and here, we obtain
$a_3 = \text{\rm Tr}  (g^* | \chi_5) = e_1$.

\par \vskip 1pc \noindent
{\bf Lemma 2.5.} (1) Suppose that $g$ switches $\chi_4$ with $\chi_4'$.
Then w.r.t. to one and the same basis $\{u_1, \dots, u_8\}$, we have the 
following matrix representation of $S_3 \times \mu_4$ on $\chi_4 \oplus \chi_4'$,
where $\eta_3$ is a primitive 3rd root of 1.
Moreover, $(g^2)^* | \chi_4 = (d_1d_5) \, \text{\rm id} = (g^2)^* | \chi_4'$:
$$\tau^* = [1, 1, \eta_3, \eta_3^2, 1, 1, \eta_3, \eta_3^2],$$
$$\sigma^* = \text{\rm diag}[\pmatrix 1 & 0 \\ 0 & -1 \endpmatrix, 
\pmatrix 0 & 1 \\ 1 & 0 \endpmatrix, \pmatrix 1 & 0 \\ 0 & -1 \endpmatrix,
\pmatrix 0 & 1 \\ 1 & 0 \endpmatrix],$$
$$g^* = \pmatrix
0 & 0 & 0 & 0 & d_5 & 0 & 0 & 0 \\
0 & 0 & 0 & 0 & 0 & d_6 & 0 & 0 \\
0 & 0 & 0 & 0 & 0 & 0 & d_7 & 0 \\
0 & 0 & 0 & 0 & 0 & 0 & 0 & d_8 \\
d_1 & 0 & 0 & 0 & 0 & 0 & 0 & 0 \\
0 & d_2 & 0 & 0 & 0 & 0 & 0 & 0 \\
0 & 0 & d_3 & 0 & 0 & 0 & 0 & 0 \\
0 & 0 & 0 & d_4 & 0 & 0 & 0 & 0 \\
\endpmatrix
.$$

\par \noindent
(2) Suppose that $g$ switches $\chi_5$ with $\chi_5'$.
Then w.r.t. to one and the same basis $\{x_1, \dots, x_{10}$ $\}$, we have the 
following matrix representation of $S_3 \times \mu_4$ on $\chi_5 \oplus \chi_5'$,
where $\eta_3$ is a primitive 3rd root of 1.
Moreover, $(g^2)^* | \chi_5 = (e_1e_6) \, \text{\rm id} = (g^2)^* | \chi_5'$:
$$\tau^* = [1, \eta_3, \eta_3^2, \eta_3, \eta_3^2, 1, \eta_3, \eta_3^2, \eta_3, \eta_3^2],$$
$$\sigma^* = \text{\rm diag}[1, \,\, \pmatrix 0 & 1 \\ 1 & 0 \endpmatrix, 
\pmatrix 0 & 1 \\ 1 & 0 \endpmatrix, 1, \,\, \pmatrix 0 & 1 \\ 1 & 0 \endpmatrix,
\pmatrix 0 & 1 \\ 1 & 0 \endpmatrix],$$
$$g^* = \pmatrix
0 & 0 & 0 & 0 & 0 & e_6 & 0   & 0   & 0  & 0\\
0 & 0 & 0 & 0 & 0 & 0   & e_7 & 0   & 0  & 0 \\
0 & 0 & 0 & 0 & 0 & 0   & 0   & e_7 & 0  & 0 \\
0 & 0 & 0 & 0 & 0 & 0   & 0   &   0 & e_9& 0 \\
0 & 0 & 0 & 0 & 0 & 0   & 0   &   0 & 0  & e_9 \\
e_1 & 0 & 0 & 0 & 0 & 0 & 0 & 0 & 0 & 0 \\
0 & e_2 & 0 & 0 & 0 & 0 & 0 & 0 & 0 & 0 \\
0 & 0 & e_2 & 0 & 0 & 0 & 0 & 0 & 0 & 0 \\
0 & 0 & 0 & e_4 & 0 & 0 & 0 & 0 & 0 & 0 \\
0 & 0 & 0 & 0 & e_4 & 0 & 0 & 0 & 0 & 0
\endpmatrix
.$$

\par \vskip 1pc \noindent
{\it Proof.} The proof is similar to {\bf (2.4)}.

\par \vskip 1pc \noindent
To prove {\bf (2.2)}, we consider first the case where
both $\chi_4$ and $\chi_5$ are $g$-stable:

\par \vskip 1pc \noindent
{\bf Lemma 2.6.} Suppose that both $\chi_4$ and $\chi_5$ are $g$-stable.

\par \noindent
(1) We have the following, where
by $\sum d_1$, etc. we mean $d_1 + d_1'$ etc :
$$\gather
\chi_{\text{\rm top}}(X^{g^{\pm}}) = 2 + \text{\rm Tr}(g^* | \chi_1 \oplus \chi_1') + 
\sum (d_1 + d_3 + e_1), \\
\chi_{\text{\rm top}}(X^{g^{-1} \tau^{\mp}}) = 
\chi_{\text{\rm top}}(X^{g \tau^{\pm}}) = 2 + \text{\rm Tr}(g^* | \chi_1 \oplus \chi_1') + 
\sum (d_1 - 2d_3 + e_1), \\
\chi_{\text{\rm top}}(X^{g^2}) = 2 +  \sum (4 d_1^2 + 5e_1^2),\\
\chi_{\text{\rm top}}(X^{g^2 \tau^{\pm}}) = 2 + \sum (d_1^2 - e_1^2).
\endgather$$

\par \noindent
(2) We have $d_1^4 = e_1^4 = (d_1')^4 = (e_1')^4 = 1$ and $d_3 \in \{\pm d_1\}$,
$d_3' \in \{\pm d_1'\}$.

\par \noindent
(3) Among six 4-th roots of 1: $e_1$, $e_1'$, $d_i$, $d_i'$ ($i = 1, 3$),
either all six of them are primitive, or exactly $e_1, e_1'$ are
primitive, or exactly the $d_i$, $d_i'$ ($i = 1, 3$) are primitive 4-th 
root of 1.

\par \noindent
(4) {\bf (2.2)} holds.

\par \vskip 1pc \noindent
{\it Proof.} (1) and (2) follow from {\bf (2.4)}.
For (3), the formula for $\chi_{\text{\rm top}}(X^{g^2})$ in (1) and its
upper bound $18$ in {\bf (1.2)}
imply that at least one of the six 4-th roots of 1 in (3) is primitive.
Now (3) is a consequence of (2) and the description of 
$\chi_{\text{\rm top}}(X^g)$ and $\chi_{\text{\rm top}}(X^{g \tau})$
in (1) which and the difference (i.e., $3 \sum d_3 = 3(d_3 + d_3')$)
of which must be real numbers (indeed, integers).

\par \vskip 1pc \noindent
To prove (4), we apply (3). If
exactly these four: $d_i, d_i'$ ($i = 1, 3$) are primitive 4-th roots of 1, 
then $\chi_{\text{\rm top}}(X^{g^2 \tau}) = 2 + (-2) - 2 < 0$,
contradicting {\bf (1.4)}.
If all these six in (3) are primitive 4-th roots of 1, then
$\chi_{\text{\rm top}}(X^g)$ and $\chi_{\text{\rm top}}(X^{g \tau})$,
given in (1) and being real numbers, must all be equal to
$2 + \text{\rm Tr}  (g^* | \chi_1 \oplus \chi_1')$; hence they are all equal to
4 -- the only possible common value of these two, by {\bf (1.4)};
but then $\chi_{\text{\rm top}}(X^{g^2 \tau}) = 2 + (-2) - (-2) = 2 < 4
= \chi_{\text{\rm top}}(X^{g \tau})$, a contradiction to {\bf (1.4)}.

\par \noindent
Thus, exactly $e_1, e_1'$ are primitive 4-th root of 1,
while $d_i, d_i' \in \{\pm 1\}$ ($i = 1, 3$).
So (*) : \,\, $-2 \le \chi_{\text{\rm top}}(X^g) \le 8$. Also
$\chi_{\text{\rm top}}(X^{g^2}) = 2 + 4 \times 2 + 5 \times (-2) = 0$ and
$\chi_{\text{\rm top}}(X^{g^2 \tau^{\pm}}) = 2 + 2  - (-2) = 6$.
Now (1) implies that 
$\chi_{\text{\rm top}}(X^{g \tau^{\pm}}) + 3 \sum d_3 = \chi_{\text{\rm top}}(X^g) = 0$
(mod 4) by {\bf (1.4)}, and also $\sum d_3 = d_3 + d_3' \in \{0, \pm 2\}$
and $\chi_{\text{\rm top}}(X^{g \tau^{\pm}}) \in \{2, 4, 6\}$ by {\bf (1.4)}.
These and the (*) above infer that the cases in {\bf (2.2)} occur. The lemma is proved.

\par \vskip 1pc \noindent
The first two assertions below are consequences of {\bf (2.4)} - {\bf (2.5)} 
and an argument similar to {\bf (2.6)}.

\par \vskip 1pc \noindent
{\bf Lemma 2.7.} Suppose that $g$ switches $\chi_4$ with $\chi_4'$
but keeps $\chi_5$ (and $\chi_5'$) stable. 

\par \noindent
(1) We have the following, where
$\delta \in S_3 = \langle \sigma, \tau \rangle$ and
by $\sum e_1$ etc. we mean $e_1 + e_1'$ etc :
$$\gather 
\chi_{\text{\rm top}}(X^{g^{-1} \delta^{-1}}) = \chi_{\text{\rm top}}(X^{g \delta}) = 
2 + \text{\rm Tr}(g^* | \chi_1 \oplus \chi_1') + \sum e_1,\\
\chi_{\text{\rm top}}(X^{g^2}) = 2 + 8d_1d_5 + 5 \sum e_1^2,\\
\chi_{\text{\rm top}}(X^{g^2 \tau^{\pm}}) = 2 + 2d_1d_5  - \sum e_1^2.
\endgather$$

\par \noindent
(2) We have $e_1^4 = (e_1')^4 = (d_1d_5)^2 = 1$.
Either $\{e_1, e_1'\} = \{\pm \sqrt{-1}\}$, or $e_1, e_1' \in \{\pm 1\}$.

\par \noindent
(3) {\bf (2.2)} holds.

\par \vskip 1pc \noindent
{\it Proof.} To prove (3), by (1) $\chi_{\text{\rm top}}(X^g)$ ($= 0$ mod 4)
and $\chi_{\text{\rm top}}(X^{g \tau})$ ($\in \{2, 4, 6\}$) are equal (see {\bf (1.4)}).
Hence they are all equal to 4.
If both $e_1, e_1'$ are in $\{\pm 1\}$, then
$\chi_{\text{\rm top}}(X^{g^2 \tau}) = 2 + 2d_1d_5 - 2 \le 2 < 4 = \chi(X^{g \tau})$,
contradicting {\bf (1.4)}.
Thus, $\{e_1, e_1'\} = \{\pm \sqrt{-1}\}$. 
By {\bf (1.4)}, we have $4 = \chi_{\text{\rm top}}(X^{g \tau}) \le
\chi_{\text{\rm top}}(X^{g^2 \tau}) = 2 + 2 d_1d_5 + 2$, whence the latter equals
$6$ and $d_1d_5 = 1$. Now $\chi_{\text{\rm top}}(X^{g^2}) = 2 + 8 + 5 \times (-2) = 0$.
Therefore, the second case in {\bf (2.2)} occurs. This proves the lemma.

\par \vskip 1pc \noindent
{\bf Lemma 2.8.} Suppose that $\chi_4$ (and $\chi_4'$ are)
is $g$-stable but
$g$ switches $\chi_5$ with $\chi_5'$.

\par \noindent
(1) We have the following, where
by $\sum d_1$ etc. we mean $d_1 + d_1'$ etc :
$$\gather
\chi_{\text{\rm top}}(X^{g^{\pm}}) = 2 + \text{\rm Tr}
(g^* | \chi_1 \oplus \chi_1') + \sum (d_1 + d_3),\\
\chi_{\text{\rm top}}(X^{g^{-1} \tau^{\mp}}) = 
\chi_{\text{\rm top}}(X^{g \tau^{\pm}}) = 2 + \text{\rm Tr}(g^* | \chi_1 \oplus \chi_1') + 
\sum (d_1 - 2d_3),\\
\chi_{\text{\rm top}}(X^{g^2}) = 2 + 4 \sum d_1^2 + 10 e_1 e_6,\\
\chi_{\text{\rm top}}(X^{g^2 \tau^{\pm}}) = 2 + \sum d_1^2 - 2 e_1 e_6.
\endgather$$

\par \noindent
(2) We have $d_1^4 = (d_1')^4 = (e_1 e_6)^2 = 1$ and $d_3 \in \{\pm d_1\}$,
$d_3' \in \{\pm d_1'\}$.

\par \noindent
(3) Either the four 4-th roots of 1: $d_i, d_i'$ ($i = 1, 3$) are
all in $\{\pm \sqrt{-1}\}$, or these four are all in $\{\pm 1\}$ (so
$e_1 e_6 = -1$ and $\chi_{\text{\rm top}}(X^{g^2}) = 0$ by {\bf (1.2)}).

\par \noindent
(4) {\bf (2.2)} holds.

\par \vskip 1pc \noindent
{\it Proof.} (1) - (2) are consequences of {\bf (2.5)} - {\bf (2.6)},
while the proof of (3) - (4) are similar to the argument for the case 
of {\bf (2.6).}
Indeed, if the first (resp. second) situation in (3) occurs,
then a contradiction (resp. {\bf (2.2)} holds). This proves the lemma.

\par \vskip 1pc \noindent
{\bf Lemma 2.9.} Suppose that $g$ switches $\chi_4$ with $\chi_4'$
and $\chi_5$ with $\chi_5'$. Then {\bf (2.2)} holds.

\par \noindent
To be precise, we have the following,
where $\delta$ is in
$S_3 = \langle \sigma, \tau \rangle$,
where $(d_1d_5)^2 = (e_1e_6)^2 = 1$:
$$\gather
\chi_{\text{\rm top}}(X^{g^{-1} \delta^{-1}}) =
\chi_{\text{\rm top}}(X^{g \delta}) = 2 + \text{\rm Tr}  (g^* | \chi_1 \oplus \chi_1'),\\
\chi_{\text{\rm top}}(X^{g^2}) = 2 + 8d_1d_5 + 10 e_1 e_6,\\
\chi_{\text{\rm top}}(X^{g^2 \tau^{\pm}}) = 2 + 2d_1d_5  - 2e_1e_6.
\endgather$$

\par \vskip 0.5pc \noindent
{\it Proof.} The formulae or equalities
are consequences of {\bf (2.4)} - {\bf (2.5)}.
As in {\bf (2.7)}, the formulae in (1) and {\bf (1.4)} imply that 
$\chi_{\text{\rm top}}(X^g) = \chi_{\text{\rm top}}(X^{g \tau}) = 4$.
The formula for $\chi_{\text{\rm top}}(X^{g^2 \tau})$ 
and its lower bounder $4 = \chi_{\text{\rm top}}(X^{g \tau})$ by {\bf (1.4)},
infer that it equals $6$ and $d_1d_5 = 1$, $e_1 e_6 = -1$.
This proves the lemma. The proof of {\bf (2.2)} is completed.

\par \vskip 2pc \noindent
{\bf \S 3. The proofs of Theorems A, B and C}

\par \vskip 1pc \noindent
In this section we shall prove Theorems A, B and C.
We first prove the result below which is a consequence
of {\bf (3.2)}-{\bf (3.8)} below.

\par \vskip 1pc \noindent
{\bf Theorem 3.1.} 

\par \noindent
(1) There is no faithful group action
of the form $A_5 . \mu_4$ (see {\bf (1.0)}) on a $K3$ surface.

\par \noindent
(2) If $A_5 . \mu_I$ acts faithfully on a $K3$ surface, then $I = 1$, or $2$.

\par \vskip 1pc \noindent
(2) follows from (1), {\bf (1.1)} and [Z2, Theorem 3.1]. Let us prove {\bf (3.1)} (1).
Suppose the contrary that $G := A_5 . \mu_4$ acts on a $K3$ surface $X$.
Then $G = A_5 : \mu_4$
and the {\it unique} group structure of such $G$ is given in {\bf (1.7)}.
We use the notation in {\bf (2.1)} and {\bf (2.2)}.
First, we need:

\par \vskip 1pc \noindent
{\bf Proposition 3.2.} Suppose that $G = A_5 : \mu_4$ acts on
a $K3$ surface $X$. Then with the notation in {\bf (2.1)}, the fixed locus
$X^{g^2} = C \coprod_{i=1}^6 D_i$ is a disjoint union of a genus-$7$ curve $C$
(hence $C^2 = 12$) and six smooth rational curves.
Both $C$ and $\sum_{i=1}^6 D_i$ are $G$-stable.

\par \vskip 1pc \noindent
{\it Proof.} We apply {\bf (2.2)}.
Then we always have $\chi_{\text{\rm top}}(X^{g^2}) = 0$.
Also {\bf (1.4)} implies that $X^{g^2} \supseteq X^{g} \ne \emptyset$,
so either $X^{g^2} = \coprod_{i=1}^s E_i$ with $1 \le s \le 10$ (by {\bf (1.2)}) is a disjoint union of a 
few smooth elliptic curves $E_i$ (so $X^g_{1-\text{\rm dim}}$ is, if not empty,
a disjoint union of some of the $E_i$'s,
and hence $n_g = 0$ in notation of {\bf (1.4)}), 
or $X^{g^2} = C \coprod_{i=1}^s D_i$ is a disjoint
union of a smooth curve $C$ and $s$ smooth rational curves $D_i$ 
with $9 \ge s = g(C) - 1 \ge 1$ (see {\bf (1.2)}).

\par \vskip 1pc \noindent
Consider the case where $X^{g^2} = \coprod_{i=1}^s E_i$.
Then $n_g = 0$ and 
$(n_g, m_g; \chi_{\text{\rm top}}(X^g), \chi_{\text{\rm top}}(X^{g \tau})$, $\chi_{\text{\rm top}}(X^{g^2 \tau})$,
$\chi_{\text{\rm top}}(X^{g^2})) = (0, 4; 4, 4, 6, 0)$. 
Note that $|X^g_{\text{\rm isol}}| =
m_g = 4$.
We may assume that $E_1$ contains an isolated $g$-fixed point.
Since the restriction $g | E_1$ is now of order 2,
this $E_1$ contains all four isolated 
$g$-fixed points by {\bf (1.9)}.
Now $g$ commutes with every element of $\langle \sigma, \tau \rangle
= S_3$ as mentioned in {\bf (2.1)}, and hence there is a natural homomorphism
$S_3 \rightarrow S_4$ ($=$ the full symmetry group of the 4-point
set $X^g_{\text{\rm isol}}$). By {\bf (1.2)} and {\bf (1.9)},
the restriction $\tau | X^g_{\text{\rm isol}} \ne \text{\rm id}$.
So the image of this homomorphism
equals one of the four 1-point (say $P_1$) stabilizer subgroups ($\cong S_3$)
in $S_4$. This leads to that $S_3 < SL(T_{X, P_1})$,
contradicting {\bf (1.0C)}.

\par \noindent
Next we consider the case where
$X^{g^2} = C \coprod_{i=1}^s D_i$.
We claim that $s = 1, 5, 6$.
Since $g^2$ is in the centre of $G$ by {\bf (2.1)},
our $G$ acts on $X^{g^2}$ and hence stabilizes $C$
and permutes $D_i$'s. Note that $C$ and the $A_5$-orbits of
$\{D_1, \dots, D_s\}$ will give linearly independent
classes in $S_X^{A_5} \otimes {\bold Q}$. Since the latter is of rank 2 by {\bf (1.6)},
this $A_5$ acts transitively on the set $\{D_1, \dots, D_s\}$
and hence the image of the natural homomorphism $A_5 \rightarrow S_s$
is a transitive subgroup of $S_s$. 
Now the claim follows from {\bf (1.8)}.

\par \noindent
We assert that $C$ is not $g$-fixed.
Indeed, let $\delta = (13)(24)$, then $c_{\delta}(g) = g \sigma$
with $\sigma = (12)(34)$ 
(because $c_g = c_{(12)}$ on $A_5$).
Hence $X^{g \sigma} = \delta (X^g)$. So $\delta(C)$ is contained in $X^{g \sigma}
\subseteq X^{g^2}$ (noting that $(g \sigma)^2 = g^2$), whence it equals
the unique curve $C$ of genus $\ge 2$ in $X^{g^2}$. Thus $C = \delta(C)$ is 
pointwise $g \sigma$-fixed. However, $C$ is also pointwise $g$-fixed,
whence is pointwise $\sigma$-fixed. This contradicts {\bf (1.2)}.
So the assertion is proved.

\par \noindent
We claim that $s = 1$ is impossible.
Consider the case $s = 1$. Then $G = A_5 : \langle g \rangle$ acts on the
set $\{C, D_1\}$ and hence stabilizes both $C$ and $D_1$.
If $D_1$ is pointwise $g$-fixed, then as above,
$D_1$ would be pointwise
($g \sigma$ and hence) $\sigma$-fixed, a contradiction.
So the restriction $g|D_1$ is not identity.
We consider the natural homomorphism 
$f : S_5 = A_5 : \langle \overline{g} \rangle = G/\langle g^2 \rangle \rightarrow \text{Aut}(D_1)$ (see {\bf (2.1)},
where $\overline{g}$ is the coset in $\langle g \rangle/\langle g^2 \rangle$
containing $g$.
Clearly, the restriction $f | A_5$ is an injection by {\bf (1.2)}.
Hence $|\text{\rm Ker}(f)| \le 2$ and $\text{\rm Ker}(f)$ is normal in $S_5$. So
$\text{\rm Ker}(f) = (1)$ and $S_5 \cong f(S_5) < \text{Aut}({\bold P}^1)$, contradicting 
{\bf (1.9)}.

\par \noindent 
We still have to rule out the case $s = 5$.
Since $C$ is not pointwise $g$-fixed as proved above,
$X^g_{1-\text{\rm dim}}$ is (if not empty) a disjoint union of $n_g/2$ ($\ge 0$) of $D_i$'s.
If $\tau = (345)$ stabilizes some $D_j$ then
$\tau$ fixes exactly two points on $D_j$ by
{\bf (1.2)} and {\bf (1.9)}. Since $|X^{\tau}| = 6$,
this $\tau$ stabilizes at most three $D_j's$. Thus we may assume that
$\tau$ permutes $D_1, D_2, D_3$ while stabilizes $D_4, D_5$.
Now the commutability of $g$ with $\tau$ implies that
$g$ stabilizes each $D_i$ ($i = 1, 2, 3$);
also none of $D_i$ ($i = 1, 2, 3$) is pointwise
$g$-fixed, for otherwise all these three $D_i$ (forming one $\tau$-orbit)
are pointwise $g$-fixed, whence $n_g \ge 3$, contradicting {\bf (2.2)}.
Thus, $m_g = |X^g_{\text{\rm isol}}| \ge |\sum_{i=1}^3 |D_i^g| = 6$.
So the first case in {\bf (2.2)} occurs and $n_g = 1$, $m_g = 6$.
Here $n_g = 1$ implies that (after switching $D_4$ with $D_5$ if necessary)
$D_5$ is poinwise $g$-fixed, and $D_4$ is $g$-stable but not
$g$-fixed. This leads to 
$6 = |X^g_{\text{\rm isol}}| \ge \sum_{i = 1}^4 |D_i^g| = 8$,
a contradiction. 
So {\bf (3.2)} is proved.
Indeed, for the last part, note that $g^2$ is in the centre of $G$ 
by {\bf (2.1)} and hence $G$ acts on $X^{g^2}$.

\par \vskip 1pc \noindent
We continue the proof of {\bf (3.1)} (1).
In notation of {\bf (3.2)},
we set $D = \sum_{i=1}^6 D_i$ and $L_0 := {\bold Z}[C, D]$.
Then we have:

\par \vskip 1pc \noindent
{\bf Lemma 3.3.} Suppose that $G = A_5 : \mu_4$ acts on a $K3$ surface $X$.

\par \noindent
(1) $L_0$ is a sublattice (with intersection form $\text{\rm diag}[12, -12]$)
of $S_X^{A_5}$ of finite index $d_1$. In particular,
$S_X^G = S_X^{A_5}$, i.e., $g^* | S_X^{A_5} = \text{\rm id}$.

\par \noindent
(2) If $d_1 > 1$, then $d_1 = 2$ and
$S_X^{A_5}$ equals ${\bold Z}[u_1, u_2]$ with $u_1 = \frac{1}{2}(C + D)$ and 
$u_2 = \frac{1}{2}(C - D)$ and
with the intersection form $U(6)$, i.e., $u_i^2 = 0$ and $u_1 . u_2 = 6$.

\par \vskip 1pc \noindent
{\it Proof.} (1) Clearly, $S_X^{A_5} \supseteq S_X^G \supseteq L_0$ 
by {\bf (3.2)}. Now (1) follows from the fact that $\text{\rm rank} \, S_X^{A_5} = 2$
by {\bf (1.6)}.

\par \noindent
(2) Suppose that $d_1 > 1$. Let $\theta = \frac{1}{12}(a C + b D)$
be in $S_X^{A_5} \subseteq L_0^{\vee} = \text{\rm Hom}(L_0, {\bold Z})
= {\bold Z}[C/12, D/12]$ but not in $L_0$. Since 
$-2b/12 = \theta . D_1 \in {\bold Z}$, we have $6 | b$.
This and $(a^2 - b^2)/12 = \theta^2 \in {\bold Z}$ imply that
$12$ divides $a^2$, whence $6 | a$. So modulo $L_0$, 
our $\theta = C/2$, or $D/2$ or $(C + D)/2$.
Since $\theta^2 \in 2{\bold Z}$, we have
$\theta = (C+D)/2$ and hence
$S_X^{A_5} = {\bold Z}[C, (C+D)/2] = {\bold Z}[(C+D)/2, (C-D)/2]$.
The lemma is proved.

\par \vskip 1pc \noindent
Set $L = H^0(X, {\bold Z})$ which contains $S_X \oplus T_X$ as a
sublattice of finite index. Also $L^{A_5}$ contains $S_X^{A_5} \oplus T_X$
as a sublattice of finite index $d$ by {\bf (1.0A-B)}.

\par \vskip 1pc \noindent
{\bf Lemma 3.4.} The quotient $L^{A_5}/(S_X^{A_5} \oplus T_X)$
is 2-elementary of order $d$ and isomorphic to $(0)$ ($d = 1$), 
${\bold Z}/(2)$ ($d = 2$) or $({\bold Z}/(2))^{\oplus 2}$ ($d = 4$).

\par \vskip 1pc \noindent
{\it Proof.} For a lattice $M$, we denote by $M^{\vee} = \text{\rm Hom}(M, {\bold Z})$
the dual and $A_M = M^{\vee}/M$ the discriminant group.
Then we have,
where $\iota$ is the inclusion:
$$\gather
S_X^{A_5} \oplus T_X \subseteq L^{A_5} \subseteq (L^{A_5})^{\vee}
\subseteq (S_X^{A_5})^{\vee} \oplus T_X^{\vee}, \\
\iota : L^{A_5}/(S_X^{A_5} \oplus T_X) \rightarrow A_{S_X^{A_5}} 
\oplus A_{T_X}.
\endgather$$
Let $pr_1$ and $pr_2$ be the projections from $A_{S_X^{A_5}} \oplus A_{T_X}$
to its first and second summands, respectively.
Since $S_X^{A_5}$ and $T_X$ are primitive in $L^{A_5}$, 
both compositions $pr_i \circ \iota$ are injective. In particular,
the quotient group in {\bf (3.4)} is regarded as a subgroup
of a bigger group $A_{T_X}$, whence is generated by 2 elements
because the same is true for the bigger group (since $\text{\rm rank} \, T_X = 2$
by {\bf (1.1)}).
We still have to show that this quotient group is 2-elementary.

\par \noindent
Take a coset $\overline{\theta}$ from the quotient group in {\bf (3.4)}.
In notation of {\bf (1.5)}, we write
$$\theta = u + \frac{1}{2m} (a t_1 + b t_2) \in 
(S_X^{A_5})^{\vee} \oplus T_X^{\vee}.$$
Regarding $\overline{\theta}$ as an element of $A_{S_X^{A_5}}$
via the injection $pr_1 \circ \iota$, we have by {\bf (3.3)}, 
modulo $S_X^{A_5} \oplus T_X$,
that 
$$0 = g^* \theta - \theta = \frac{1}{2m}[a(g^*t_1 - t_1) + b(g^*t_2 - t_2)]
= \frac{1}{2m} [-(a+b) t_1 + (a-b) t_2].$$
So $2m$ divides $a + b$, $a - b$ (and hence $2a$ and $2b$) because $T_X$
is primitive in $L$. Thus $m$ divides $a$ and $b$ and we write
$a = m a'$ and $b = m b'$ so that
$\theta = u + \frac{1}{2}(a' t_1 + b' t_1)$.
Therefore, modulo $T_X$, we have $2 u = 2 \theta \in 2 L^{G_N} \subset L^{G_N}$,
whence $2 u \in L \cap (S_X^{A_5})^{\vee} = S_X^{A_5}$ (because
the latter is primitive in $L$). So $2 \overline{\theta} = 0$.
The lemma is proved.

\par \vskip 1pc \noindent
{\bf Lemma 3.5.} One of the following cases occurs.

\par \noindent
(1) We have $m = 5$. Both the quotients $S_X^{A_5}/L_0$ and
$L^{A_5}/(S_X^{A_5} \oplus T_X)$ are isomorphic to ${\bold Z}/(2)$.
Moreover, the discriminant form of $(L^{A_5})^{\vee}/L^{A_5}
\cong ({\bold Z}/(30))^{\oplus 2}$ is given in [Z2, Theorem 2.1]
(corresponding to the matrix $M_1$ there) and generated by the cosets
$\overline{\varepsilon}_i$ with
$\varepsilon_1 = e_1^*, \, \varepsilon_2 = e_2^* + e_3^* + e_4^*$
and the intersection form
(note that $\overline{\varepsilon}_i^2$ is in ${\bold Q}/2{\bold Z}$
while $\overline{\varepsilon}_1 . \overline{\varepsilon}_2$ is in ${\bold Q}/{\bold Z}$):
$$(\overline{\varepsilon}_i . \overline{\varepsilon}_j) = \pmatrix -23/30 & -1/5 \\ -1/5 & -35/30 \endpmatrix.$$

\par \noindent
(2) We have $m = 10$, $S_X^{A_5}/L_0 \cong {\bold Z}/(2)$ and
$L^{A_5}/(S_X^{A_5} \oplus T_X) \cong {(\bold Z}/(2))^{\oplus 2}$.

\par \noindent
(3) We have $m = 5$, $L_0 = S_X^{A_5}$ and
$L^{A_5}/(S_X^{A_5} \oplus T_X) \cong {(\bold Z}/(2))^{\oplus 2}$.

\par \vskip 1pc \noindent
{\it Proof.} In notation of {\bf (3.3)} and {\bf (3.4)}, 
we have $-(12^2)(4m^2) = |L_0||T_X| = d_1^2 d^2 |L^{A_5}|$.
On the other hand, 
$-|L^{A_5}| = 30^2, 3 \times 10^2, 20^2, 3 \times 20^2, 3 \times 40^2$
by the calculation in [Z2, Theorem 2.1]. Then the lemma follows easily.

\par \vskip 1pc \noindent
{\bf Lemma 3.6.} The case (3) in {\bf (3.5)} does not occur.

\par \vskip 1pc \noindent
{\it Proof.} Consider the case (3) in {\bf (3.5)}. 
Let $\theta$ be an element in $L^{A_5}$ but not in the smaller set
$S_X^{A_5} \oplus T_X$.
We claim that $\theta^2 \in 2{\bold Z}$ implies that 
modulo this smaller set, our $\theta$ equals some $\theta_j$ below,
where $u_1 := C$, $u_2 := D$ and $T_X = {\bold Z}[t_1, t_2]$ as
in {\bf (1.5)}. Here
$\theta_j := \frac{1}{2}(t_1 + t_2) + \frac{1}{2}u_j$.

\par \noindent
Indeed, since the quotient group in {\bf (3.5) (3)} is 2-elementary,
we can write, modulo the smaller set, 
that $\theta = \frac{1}{2}(a_1 t_1 + a_2 t_2 + b_1 u_1 + b_2 u_2)$
with $a_i$, $b_j$ in $\{0, 1\}$ but not all zero.
Indeed, $(a_1, a_2) \ne (0, 0) \ne (b_1, b_2)$ because
both $S_X^{A_5}$ and $T_X$ are primitive in $L$.
Now modulo $2{\bold Z}$, we have the following,
so the claim follows:
$$\frac{1}{2}(a_1^2 + a_2^2) + b_1^2 + b_2^2 = \frac{2m}{4}(a_1^2 + a_2^2)
+ \frac{12}{4} (b_1^2 - b_2^2) = \theta^2 = 0.$$

\par \noindent
Since $\theta_1 - \theta_2$ is not in $L^{A_5}$
(not in $L$ at all, by the primitivity of $S_X^{A_5}$ in $L$),
at most one of $\theta_j$ is in $L^{A_5}$. So $L^{A_5}/(S_X^{A_5} \oplus T_X)$
is of order $\le 2$, a contradiction.

\par \vskip 1pc \noindent
We start anew. By {\bf (3.3)} and {\bf (3.6)}, the lattice
$S_X^{A_5}$ equals ${\bold Z}[u_1, u_2]$ with $u_1 = \frac{1}{2}(C + D)$
and $u_2 = \frac{1}{2}(C - D)$, and has the intersection form $U(6)$.

\par \vskip 1pc \noindent
{\bf Lemma 3.7.} The case (2) in {\bf (3.5)} is impossible.

\par \vskip 1pc \noindent
{\it Proof.} Take $\theta$ in $L^{A_5}$ but not in the smaller set
$S_X^{A_5} \oplus T_X$. As in {\bf (3.6)}, $\theta^2 \in 2{\bold Z}$ implies that
modulo the smaller set, our $\theta$ is one of the following
$$\theta^i = \frac{1}{2} t_i + \frac{1}{2}(u_1 + u_2), \hskip 1pc
\theta_j = \frac{1}{2}(t_1 + t_2) + \frac{1}{2} u_j.$$

\par \noindent
Since $\theta^1 - \theta^2$ is not in $L^{A_5}$ (not in $L$ at all),
not both $\theta^i$ are in $L^{A_5}$. By the same reasoning
not both $\theta_j$ are in $L^{A_5}$.
Since $L^{A_5}/(S_X^{A_5} \oplus T_X) \cong ({\bold Z}/(2))^{\oplus 2}$
is generated by two elements,
one of $\theta^i$ ($i = 1, 2$) and one  of 
$\theta_j$ ($j = 1, 2$) are in $L^{A_5}$. But
$\theta^i . \theta_j = \frac{2m}{4} + \frac{6}{4} = \frac{13}{2}$,
which is not an integer. This is a contradiction.

\par \vskip 1pc \noindent
{\bf Lemma 3.8.} Suppose the case (1) in {\bf (3.5)} occurs. Then
we have:

\par \noindent
(1) $L^{A_5}$ is generated by $S_X, T_X$ and 
$\theta = \frac{1}{2}(t_1 + t_2 + u_1 + u_2)$.

\par \noindent
(2) The discriminant group $A_{L^{A_5}} = (L^{A_5})^{\vee}/L^{A_5}$
(with the dual $(L^{A_5})^{\vee} = \text{\rm Hom}(L^{A_5}, {\bold Z})$)
is generated by the cosets $\overline{\delta}_j$ ($j = 1, 2$) 
which (together with the intersection form) is given as follows:
(where $t_i^* . t_j = \delta_{ij}$, and $u_i^* . u_j = \delta_{ij}$
in Kronecker's symbol):
$$\delta_1 = t_2^* + u_1^* + 2 u_2^* = 
\frac{1}{10} t_2 + \frac{1}{6} (2u_1 + u_2), \hskip 1pc
\delta_2 = t_1^* + u_1^* = \frac{1}{10} t_1 + \frac{1}{6} u_2,$$
$$(\overline{\delta}_i . \overline{\delta}_j) = \pmatrix
23/30 & 1/3 \\ 1/3 & 1/10 \endpmatrix.$$

\par \vskip 1pc \noindent
{\it Proof.} (1) can be proved as in {\bf (3.6)}, by making use
of that $\theta_1^2 \in 2{\bold Z}$ for every $\theta_1$
in $L^{A_5}$.

\par \noindent
(2) Since $\delta_i . \theta$, $\delta_i . t_j$ and $\delta_i . u_j$ are
all in ${\bold Z}$ by a direct calculation, we see that both $\delta_i$
are in $(L^{A_5})^{\vee}$. One checks easily that the subgroup $\langle \overline{\delta}_1, 
\overline{\delta}_2 \rangle$ of the discriminant group in (2) is isomorphic to
$({\bold Z}/(30))^{\oplus 2}$, whence this subgroup is indeed the whole
discriminant group in (2) (because the latter is of order $30^2$ by {\bf (3.5)}).
This proves the lemma.

\par \vskip 1pc \noindent
Here comes the punch line. By {\bf (3.5)} - {\bf (3.8)},
there is an isometry $\varphi : \langle \overline{\varepsilon}_1, 
\overline{\varepsilon}_2 \rangle
\longrightarrow \langle \overline{\delta}_1, \overline{\delta}_2 \rangle$,
so for some integers $a, b, c, d$,
we can write $(\varphi(\overline{\varepsilon}_1), \varphi(\overline{\varepsilon}_2))
= (\overline{\delta}_1, \overline{\delta}_2) \pmatrix a & c \\ b & d \endpmatrix$.
Thus,
$$\gather
-23/30 = \varepsilon_1^2 = \varphi(\varepsilon_1)^2 =
(a \delta_1 + b \delta_2)^2 = \frac{1}{30}(23 a^2 + 3 b^2 + 20ab) \,\, 
(\text{\rm mod} \,\, 2{\bold Z}), \\
-23 = 23 a^2 + 3 b^2 + 20ab \,\, (\text{\rm mod} \,\, 60{\bold Z}).
\endgather$$

\par \noindent 
The congruence above implies that modulo 4, we have
$1 = -a^2 - b^2$, which is impossible.
This completes the proof of {\bf (3.1)} (1) and also the whole of {\bf (3.1)}.

\par \vskip 1pc \noindent
We now prove Theorems A, B and C in the Introduction.
In Theorem C, we have $H \le G_N$ by {\bf (1.0A)}; so $H$
is either one of $A_5$, $L_2(7)$,
$A_6$ and $M_{20} = C_2^{\oplus 4} : A_5$, by [Xi, the list];
if $H = L_2(7)$ then $G_N = H$ by [Mu1] and Theorem C follows from
[OZ3, Main Theorem].

\par \noindent
Therefore, we may assume that in all three theorems,
$G$ is a finite group containing $A_5$ and acting faithfully on
a $K3$ surface $X$.
Write $G = G_N . \mu_I$ as in {\bf (1.0)}. By {\bf (1.0A)} the
$A_5$ in $G$ is contained in $G_N$.
So $G_N$ is either one of $A_5$, $S_5$,
$A_6$ and $M_{20} = C_2^{\oplus 4} : A_5$, by [Xi, the list]. 

\par \noindent
Consider the case $G_N = A_5$. Then $I = 1, 2$,
by {\bf (1.1)}, [Z2, Theorem 3.1] and {\bf (3.1)}. If $I = 1$, then $G = A_5$.
If $I = 2$, let $\rho : G \rightarrow S_5 \times \mu_2$
($x \mapsto (c_x, \alpha(x))$) be the injection as in {\bf (1.7)}
so that $\rho(A_5) = A_5 \times \langle 1 \rangle$;
if the projection $pr_1 : S_5 \times \mu_2 \rightarrow S_5$ maps
$\rho(G)$ to $A_5$ (resp. to $S_5$), then $G \cong \rho(G) = A_5 \times \mu_2$
(resp. $G \cong \rho(G) \cong pr_1(\rho(G)) = S_5$, by comparing the orders);
see the argument below.
Thus Theorems A, B and C are true.

\par \noindent
Consider the case $G_N = S_5$. Let $g$ be in $G$ such that $\alpha(g)$ is 
a generator of $\mu_I$. Since $\text{Aut}(S_5) = S_5$ and
$x \mapsto c_x$ gives rise to an isomorphism $S_5 \rightarrow \text{Aut}(S_5)$, 
we see that the map $G \rightarrow \text{Aut}(S_5) = S_5$ ($x \mapsto c_x$)
is surjective, and the conjugation maps $c_g = c_s$ on $S_5$, for some $s \in S_5$.
Replacing $g$ by $g s^{-1}$, we may assume that $g$ commutes with
every element in $G_N = S_5$. So $g^I \in \text{\rm Ker}(\alpha) = G_N$ is in the centre of $G_N = S_5$ 
(which is $(1)$), whence $\text{\rm ord} (g) = I$, while $\alpha(g)$ is a generator of
$\mu_I$.
Thus $G = S_5 \times \mu_I \ge A_5 \times \mu_I$.
So $I = 1, 2$ by {\bf (1.1)}, [Z2, Theorem 3.1] and {\bf (3.1)}.
Hence Theorems A, B and C are true.

\par \noindent
Consider the case where $G_N = A_6$ or $G_N = M_{20} = C_2^4 : A_5$.
Then $G_N$ does not contain an $A_5$ as a normal subgroup
(otherwise, in the latter case, $M_{20} = C_2^4 \times A_5$, absurd).
So $A_5$ is also not normal in $G$. Thus Theorems A and B are void this time.
Now Theorem C follows from [Ko2] and [KOZ1].

\par \vskip 2pc \noindent
{\bf References}

\par \noindent
[AS2] M. F. Atiyah and G. B. Segal,
The index of elliptic operators. II. 
Ann. of Math. 87 (1968), 531--545.

\par \noindent
[AS3] M. F. Atiyah and I. M. Singer,
The index of elliptic operators. III. 
Ann. of Math. 87 (1968) 546--604.

\par \noindent
[Atlas] J. H. Conway, R. T. Curtis, S. P. Norton, R. A. Parker and R. A. Wilson, 
Atlas of finite groups. Oxford University Press. Reprinted 2003 (with corrections). 

\par \noindent
[CS] J. H. Conway and N. J. A. Sloane, Sphere packings, lattices and groups. 3rd ed. 
Grundlehren der Mathematischen Wissenschaften, 290. Springer-Verlag, New York, 1999.

\par \noindent
[EDM] Encyclopedic dictionary of mathematics. Vol. I--IV. Translated from the Japanese. 
2nd ed. Edited by Kiyosi Itô. MIT Press, Cambridge, MA, 1987.

\par \noindent
[IOZ] A. Ivanov, K. Oguiso and D. -Q. Zhang, The monster and $K3$ surfaces, in preparation.

\par \noindent
[KOZ1] J. Keum, K. Oguiso and D. -Q. Zhang, The alternating group of degree 6
in geometry of the Leech lattice and $K3$ surfaces, 
Proc. London Math. Soc. 90 (2005), 371 - 394.

\par \noindent
[KOZ2] J. Keum, K. Oguiso and D. -Q. Zhang, Extensions of 
the alternating group of degree 6 in geometry of $K3$ surfaces,
math.AG/0408105, European J. Combinatorics: Special issue on Groups and
Geometries, to appear.

\par \noindent
[Ko1] S. Kondo, Niemeier lattices, Mathieu groups, and finite groups of symplectic 
automorphisms of $K3$ surfaces. 
Duke Math. J. 92 (1998), 593--598.

\par \noindent
[Ko2] S. Kondo, The maximum order of finite groups of automorphisms of $K3$ surfaces.
Amer. J. Math. 121 (1999), 1245--1252.

\par \noindent
[MO] N. Machida and K. Oguiso,
On $K3$ surfaces admitting finite non-symplectic group actions.  
J. Math. Sci. Univ. Tokyo 5 (1998), 273--297.

\par \noindent
[Mu1] S. Mukai, Finite groups of automorphisms of $K3$ surfaces and the Mathieu group. 
Invent. Math. 94 (1988), 183--221.

\par \noindent
[Mu2] Lattice-theoretic construction of symplectic actions on $K3$ surfaces,
Duke Math. J. 92 (1998), 599--603. As the Appendix to [Ko1].

\par \noindent
[Ni1] V. V. Nikulin, Finite automorphism groups of Kahler $K3$ surfaces,
Trans. Moscow Math. Soc. 38 (1980), 71--135.

\par \noindent
[Ni2] V. V. Nikulin, Factor groups of groups of automorphisms of hyperbolic forms
with respect to subgroups generated by 2-reflections. Algebrogeometric applications.
J. Soviet Math. 22 (1983), 1401--1475.

\par \noindent
[Ni3] V. V. Nikulin,
Integer symmetric bilinear forms and some of their applications. 
Math. USSR Izvestija. 14 (1980), 103 -- 167.

\par \noindent
[Og] K. Oguiso, A characterization of the Fermat quartic $K3$ surface by
means of finite symmetries, math.AG/0308062. Compositio Math. to appear.

\par \noindent
[OZ1] K. Oguiso and D. -Q. Zhang, 
On the most algebraic $K3$ surfaces and the most extremal log Enriques surfaces. Amer. J. Math. 118 (1996), 1277--1297.

\par \noindent
[OZ2] K. Oguiso and D. -Q. Zhang,
On Vorontsov's theorem on $K3$ surfaces with non-symplectic group actions. 
Proc. Amer. Math. Soc. 128 (2000), 1571--1580.

\par \noindent
[OZ3] K. Oguiso and D. -Q. Zhang, 
The simple group of order 168 and $K3$ surfaces.
Complex geometry (Gottingen, 2000), Collection of papers
dedicated to Hans Grauert, 165--184, 
Springer, Berlin, 2002. 

\par \noindent
[Sh1] I. Shimada, Rational double points on supersingular $K3$ surfaces,
Mathematics of Computation, 73 (2004), 1989--2017.

\par \noindent
[Sh2] I. Shimada, Lattices of algebraic cycles on Fermat varieties in positive characteristics,
Proc. London Math. Soc. 82 (2001), 131--172. 

\par \noindent
[Sh3] I. Shimada, On elliptic $K3$ surfaces, Michigan Math. J. 47 (2000), 423--446. 

\par \noindent
[Su] M. Suzuki, Group theory. I. 
Translated from the Japanese by the author. 
Grundlehren der Mathematischen Wissenschaften 247. 
Springer-Verlag, Berlin-New York, 1982.

\par \noindent
[Vi] E. B. Vinberg,
The two most algebraic $K3$ surfaces. 
Math. Ann. 265 (1983), 1--21.

\par \noindent
[Xi] G. Xiao, Galois covers between $K3$ surfaces.
Ann. Inst. Fourier (Grenoble) 46 (1996), 73--88.

\par \noindent
[Z1] D. -Q. Zhang, Quotients of $K3$ surfaces modulo involutions. 
Japan. J. Math. (N.S.) 24 (1998), 335--366.

\par \noindent
[Z2] D. -Q. Zhang, Niemeier lattices and K3 groups, math.AG/0408106,
Proc. Intern. Conf. Alg. Geom. in honour of Prof. Dolgachev,
Contemporary Math. Amer. Math. Soc. J. Keum (ed) to appear.

\par \vskip 2pc \noindent
D. -Q. Zhang
\par \noindent
Department of Mathematics
\par \noindent
National University of Singapore
\par \noindent
Singapore
\par \noindent
E-mail : MATZDQ$\@$MATH.NUS.EDU.SG

\enddocument